\nonstopmode
\documentclass[%
10pt]{amsart}
\usepackage{latexsym}
\usepackage{fancyhdr}
\usepackage{amsmath, amssymb}
\usepackage[ansinew]{inputenc}
\usepackage[all]{xy}
\newdir{ >}{{}*!/-10pt/@{>}}%
\newdir^{ (}{{}*!/-5pt/@^{(}}%
\newdir_{ (}{{}*!/-5pt/@_{(}}%
\usepackage{verbatim}%
\usepackage{psfrag}
\usepackage{epsfig}
\usepackage{hyperref}
\usepackage{pdfsync}
\theoremstyle{plain}
\newtheorem*{theorem*}{Theorem}
\newtheorem{theorem}{Theorem}[section]
\newtheorem*{lemma*}{Lemma}
\newtheorem{lemma}[theorem]{Lemma}
\newtheorem*{proposition*}{Proposition}
\newtheorem{proposition}[theorem]{Proposition}
\newtheorem*{corollary*}{Corollary}
\newtheorem{corollary}[theorem]{Corollary}

\theoremstyle{definition}
\newtheorem*{hypothesis*}{Hypothesis}

\newtheorem*{definition*}{Definition}
\newtheorem{definition}[theorem]{Definition}
\newtheorem*{example*}{Example}
\newtheorem{example}[theorem]{Example}
\newtheorem*{remark*}{Remark}
\newtheorem*{remarks*}{Remarks}
\newtheorem{remark}[theorem]{Remark}
\newtheorem{remarks}[theorem]{Remarks}

\newenvironment{demo}[1]{\par\smallskip\noindent{\bf #1.}}{\par\smallskip}

\allowdisplaybreaks

\newcommand{\East}[2]{-\raisebox{0.1pt}{$\mkern-16mu\frac{\;\;#1\;}{\;\;#2\;}\mkern-16mu$}\to}

\let\on\operatorname

\newcommand{\sr}[1]%
{\ifmmode{}^\dagger\else${}^\dagger$\fi\ifvmode
\vbox to 0pt{\vss
\hbox to 0pt{\hskip\hsize\hskip1em
\vbox{\hsize3cm\raggedright\pretolerance10000
\noindent #1\hfill}\hss}\vss}\else
\vadjust{\vbox to0pt{\vss%
 \hbox to 0pt{\hskip\hsize\hskip1em%
 \vbox{\hsize3cm\raggedright\pretolerance10000%
 \noindent #1\hfill}\hss}\vss}}\fi%
}

\usepackage{xstring}
\IfEndWith*{\jobname}{BEAMER}{

}{

}
\let\iold\i
\renewcommand{\i}{\ifmmode^{-1}\else\iold\fi}

\ifcase\pdfoutput

\else
\hypersetup{
bookmarksnumbered=true,%
bookmarksopen=true,%
bookmarksopenlevel=1
}
\fi

\renewcommand{\o}{\operatorname{\circ}}

\providecommand{\mapsfrom}{\kern.2em%
\setbox0=\hbox{$\leftarrow$\kern-.10em\rule[0.26mm]{0.1mm}{1.3mm}}\box0%
\kern.3em}

\allowdisplaybreaks
\title[Denjoy--Carleman Mappings of Beurling and Roumieu Type]
{The Convenient Setting for\\ Denjoy--Carleman Differentiable Mappings\\ of Beurling and Roumieu Type}

\author
[A.~Kriegl, P.W.~Michor, A.~Rainer]
{Andreas Kriegl, Peter W. Michor, and Armin Rainer}

\address{Andreas Kriegl: Fakult\"at f\"ur Mathematik, Universit\"at Wien, 
Oskar-Morgenstern-Platz~1, A-1090 Wien, Austria}
\email{andreas.kriegl@univie.ac.at}

\address{Peter W. Michor: Fakult\"at f\"ur Mathematik, Universit\"at Wien, 
Oskar-Morgenstern-Platz~1, A-1090 Wien, Austria}
\email{peter.michor@univie.ac.at}

\address{Armin Rainer: Fakult\"at f\"ur Mathematik, Universit\"at Wien, 
Oskar-Morgenstern-Platz~1, A-1090 Wien, Austria}
\email{armin.rainer@univie.ac.at}

\thanks{AK was supported by FWF-Project P~23028-N13; PM by FWF-Project P~21030-N13; 
AR by FWF-Project P~22218-N13}

\begin{document}

\begin{abstract}
We prove in a uniform way that all Denjoy--Carleman differentiable function classes of Beurling type $C^{(M)}$ 
and of Roumieu type $C^{\{M\}}$, 
admit a convenient setting if 
the weight sequence $M=(M_k)$ is log-convex and of moderate growth: 
For $\mathcal C$ denoting either $C^{(M)}$ or $C^{\{M\}}$, 
the category of $\mathcal C$-mappings is cartesian closed in the sense that 
$\mathcal C(E,\mathcal C(F,G))\cong \mathcal C(E\times F, G)$ for convenient vector spaces.
Applications to manifolds of mappings are given: The group of
$\mathcal C$-diffeomorphisms is a regular $\mathcal C$-Lie group if $\mathcal C \supseteq C^\omega$, but not better. 
\end{abstract}

\maketitle

\section{\label{nmb:1}Introduction}

Denjoy--Carleman differentiable functions form classes of smooth functions that are 
described by growth conditions on the Taylor expansion. 
The growth is prescribed in terms of a sequence $M=(M_k)$ of positive real numbers which 
serves as a weight for the iterated derivatives: for compact $K$ the sets 
\[
\Big\{\frac{f^{(k)}(x)}{\rho^{k} \, k! \, M_{k}} : x \in K, k \in \mathbb{N} \Big\}
\] 
are required to be bounded. The positive real number $\rho$ is subject to either 
a universal or an existential quantifier, thereby dividing the Denjoy--Carleman classes into 
those of Beurling type, denoted by $C^{(M)}$, and those of Roumieu type, denoted by $C^{\{M\}}$, 
respectively. For the constant sequence $M=(M_k)=(1)$, as Beurling type we recover the real and imaginary parts 
of all entire functions on the one hand, and as Roumieu type the real analytic functions on the other hand, 
where $1/\rho$ plays the role of a radius of convergence.
Moreover, Denjoy--Carleman classes are divided into quasianalytic and non-quasianalytic classes,
depending on whether the mapping to infinite Taylor expansions is injective on the class or not.

That a class of mappings $\mathcal C$ admits a convenient setting means essentially that we can extend the 
class to mappings between admissible infinite dimensional spaces $E,F,\dots$ so that 
$\mathcal C(E,F)$ is again admissible 
and we have $\mathcal C(E\times F, G)$ canonically $\mathcal C$-diffeomorphic to 
$\mathcal C(E,\mathcal C(F,G))$. This property is called the \emph{exponential law}; it includes 
the basic assumption of variational calculus. 
Usually the exponential law comes hand in hand with (partially nonlinear) uniform boundedness 
theorems which are easy $\mathcal C$-detection principles.

The class $C^\infty$ of smooth mappings admits a convenient setting. This is due originally to 
\cite{Froelicher80}, \cite{Froelicher81}, and \cite{Kriegl82}, \cite{Kriegl83}.
For the $C^\infty$ convenient setting one can test smoothness along smooth curves.
Also real analytic ($C^\omega$) mappings admit a convenient setting, due to \cite{KrieglMichor90}: A mapping is $C^\omega$ if and only if it is 
$C^\infty$ and in addition is weakly $C^\omega$ along weakly $C^\omega$-curves (i.e., curves whose compositions 
with any bounded linear functional are $C^\omega$); indeed, it suffices to test along affine lines 
instead of weakly $C^\omega$-curves. 
See the book \cite{KM97} for a comprehensive treatment, or the three appendices in \cite{KMRc} for 
a short overview of the $C^\infty$ and $C^\omega$ cases.
We shall use convenient calculus of $C^\infty$-mappings in this paper, and we shall reprove that 
$C^\omega$ admits a convenient calculus. 

We now describe what was known about convenient settings for Denjoy-Carleman classes before: 
In \cite{KMRc} we developed the convenient setting for non-quasianalytic log-convex 
Denjoy--Carleman classes of Roumieu type $C^{\{M\}}$ having moderate growth, and we showed that 
moderate growth and a condition that guarantees stability under composition (like log-convexity) 
are necessary. There a mapping is $C^{\{M\}}$ if and only if
it is weakly $C^{\{M\}}$ along all weakly $C^{\{M\}}$-curves. The method of proof relies on the existence of 
$C^{\{M\}}$ partitions of unity.

We succeeded in \cite{KMRq} to prove that \emph{some} quasianalytic log-convex Denjoy--Carleman classes of 
Roumieu type $C^{\{M\}}$ having moderate growth admit a convenient setting.
The method consisted of representing $C^{\{M\}}$ as the intersection of all larger non-quasianalytic 
log-convex classes $C^{\{L\}}$. A mapping is $C^{\{M\}}$ if and only if it is weakly $C^{\{L\}}$ 
along each weakly $C^{\{L\}}$-curve for each non-quasianalytic log-convex $L \ge M$. 
We constructed countably many
classes which satisfy all these requirements, but many reasonable quasianalytic classes $C^{\{M\}}$, 
like the real analytic class, are not covered by this approach. 

In this paper we prove that \emph{all} log-convex Denjoy--Carleman classes of moderate growth 
admit a convenient setting. This is achieved through a change of philosophy: instead of \emph{testing along 
curves} as in our previous approaches \cite{KMRc} and \cite{KMRq} we \emph{test along Banach plots}, 
i.e., mappings of the respective weak class defined in open subsets of Banach spaces. By `weak' we 
mean: the mapping is in the class after composing it with any bounded linear functional.
In this way we are able to treat all Denjoy--Carleman classes uniformly, no matter if quasianalytic, 
non-quasianalytic, of Beurling, or of Roumieu type, including $C^\omega$ and real and imaginary parts 
of entire functions. Furthermore, it makes the proofs shorter and more transparent.

Smooth mappings between Banach spaces are $C^{(M)}$ or $C^{\{M\}}$ if their derivatives satisfy the 
boundedness conditions alluded to above. A smooth mapping between admissible locally convex vector spaces is 
$C^{(M)}$ or $C^{\{M\}}$ if and only if it maps Banach plots of the respective class to Banach plots 
of the same class. This implies stability under composition, see Theorem \ref{nmb:4.11}.

We equip the spaces of $C^{(M)}$ or $C^{\{M\}}$ mappings between Banach spaces with natural 
locally convex topologies which are just the usual ones if the involved Banach spaces are finite dimensional, see Section \ref{nmb:4.1}.
In order to show completeness we need to work with Whitney jets on compact subsets of Banach spaces satisfying 
growth conditions of Denjoy--Carleman type, see Proposition \ref{nmb:4.2}. 
Having found nothing in the literature we introduce Whitney jets on Banach spaces in Section \ref{nmb:3}. 

In Theorem \ref{nmb:7.2} we show that the Roumieu type classes of Denjoy--Carleman differentiable 
mappings studied in the 
present paper coincide bornologically with the classes considered previously in \cite{KMRc} and \cite{KMRq} 
and, most notably, with the structure $C^\omega$ of real analytic 
mappings introduced in \cite{KrieglMichor90} (see also \cite{KM97}). We want to stress that 
thereby we provide a considerably simpler proof for the real analytic convenient setting.
But for the results that testing along curves suffices one still has to rely on 
\cite{KrieglMichor90}, \cite{KMRc}, and \cite{KMRq}. 

For a class of mappings $\mathcal C$ that admits a convenient setting 
one can hope that the space $\mathcal C(A,B)$ of all $\mathcal C$-mappings between 
finite dimensional $\mathcal C$-manifolds (with $A$ compact for simplicity) is again a $\mathcal C$-manifold, that 
composition is $\mathcal C$, and that the group $\on{Diff}^{\mathcal C}(A)$ of all $\mathcal C$-diffeomorphisms of 
$A$ is a regular infinite dimensional $\mathcal C$-Lie group. 
In Section \ref{nmb:9} this is proved for all log-convex Denjoy--Carleman classes of moderate growth $C^{\{M\}}$ 
and for the classes $C^{(M)}$ containing $C^\omega$.

A further area of application is the perturbation theory for linear unbounded operators; see 
\cite{KMRp} and \cite{RainerN}.

This paper is organized as follows.
In Section \ref{nmb:2} we recall basic facts about Denjoy--Carleman classes $C^{[M]}$ (which stands for $C^{\{M\}}$ or $C^{(M)}$) 
in finite dimensions and discuss corresponding 
sequence spaces. 
In Section \ref{nmb:3} we introduce Whitney jets on Banach spaces. 
In Section \ref{nmb:4} we define $C^{[M]}$-mappings in infinite dimensions, first between Banach spaces with the aid of 
jets and then between convenient vector spaces, and we show that they form a category, if $M=(M_k)$ is log-convex. 
In Section \ref{nmb:5} we prove that this category is cartesian closed, if $M=(M_k)$ has moderate growth.
In Section \ref{nmb:6} we show the $C^{[M]}$ uniform boundedness principle.
In Section \ref{nmb:7} we prove that the structures studied in this paper coincide bornologically with the structures 
considered in our previous work \cite{KMRc}, \cite{KMRq}, \cite{KrieglMichor90}, and \cite{KM97}.
In Section \ref{nmb:8} we further study the spaces of $C^{[M]}$-mappings.
In Section \ref{nmb:9} we apply this theory to prove that the space of $C^{[M]}$-mappings between finite dimensional (compact) manifolds 
is naturally an infinite dimensional $C^{[M]}$-manifold, and that the group of $C^{[M]}$-diffeomorphisms of a compact manifold is a $C^{[M]}$-regular 
Lie group.   

\subsection*{Notation}
We use $\mathbb{N} = \mathbb{N}_{>0} \cup \{0\}$.
For each multi-index $\alpha=(\alpha_1,\ldots,\alpha_n) \in \mathbb{N}^n$, we write
$\alpha!=\alpha_1! \cdots \alpha_n!$, $|\alpha|= \alpha_1 +\cdots+ \alpha_n$, and 
$\partial^\alpha=\partial^{|\alpha|}/\partial x_1^{\alpha_1} \cdots \partial x_n^{\alpha_n}$. 

A sequence $r=(r_k)$ of reals is called increasing if $r_k \le r_{k+1}$ for all $k$.

We write $f^{(k)}(x) = d^k f(x)$ for the $k$-th order Fr\'echet derivative of $f$ at $x$;
by $d_v^k$ we mean $k$ times iterated directional derivatives in direction $v$.

For a convenient vector space $E$ and a closed absolutely convex bounded subset $B \subseteq E$, we denote by $E_B$ the 
linear span of $B$ equipped with the Minkowski functional $\|x\|_B = \inf \{\lambda>0 : x \in \lambda B\}$. Then $E_B$ is a 
Banach space. If $U \subseteq E$ then $U_B := i_B^{-1}(U)$, where $i_B : E_B \to E$ is the inclusion 
of $E_B$ in $E$. 

We denote by $E^*$ (resp.\ $E'$) the dual space of continuous (resp.\ bounded) linear functionals. $L(E_1,\ldots,E_k;F)$
is the space of $k$-linear bounded mappings $E_1 \times \cdots \times E_k \to F$; if $E_i =E$ for all $i$, we also write $L^k(E,F)$. 
If $E$ and $F$ are Banach spaces, then $\|~\|_{L^k(E,F)}$ denotes the operator norm on $L^k(E,F)$.
By $L_{\text{sym}}^k(E,F)$ we denote the subspace of symmetric $k$-linear bounded mappings.
We write $oE$ for the open unit ball in a Banach space $E$. 

The notation $C^{[M]}$ stands \emph{locally constantly} for either $C^{(M)}$ or $C^{\{M\}}$; this means: 
Statements that involve more than one $C^{[M]}$ symbol must not be interpreted by mixing $C^{(M)}$ and $C^{\{M\}}$.

From Section \ref{nmb:2.6} on, if not specified otherwise, a positive sequence $M=(M_k)$ is assumed to satisfy $M_0=1\le M_1$. 
In Section \ref{nmb:9} we also assume that $M=(M_k)$ is log-convex and has moderate growth, 
and in the Beurling case $C^{[M]}=C^{(M)}$ we additionally require $C^{\omega} \subseteq C^{(M)}$.

\section{\label{nmb:2}Denjoy--Carleman differentiable functions in finite dimensions}

\subsection{\label{nmb:2.1}Denjoy--Carleman differentiable functions of Beurling and Roumieu type in finite dimensions}

Let $M=(M_k)_{k \in \mathbb{N}}$ be a sequence 
of positive real numbers.
Let $U \subseteq \mathbb{R}^n$ be open, $K \subseteq U$ compact, and $\rho>0$.
Consider the set 
\begin{equation} \label{eq:1}
\Big\{\frac{\partial^\alpha f(x)}{\rho^{|\alpha|} \, |\alpha|! \, M_{|\alpha|}} : x \in K, \alpha \in \mathbb{N}^n \Big\}.
\end{equation}
We define the \emph{Denjoy--Carleman classes} 
\begin{align*}
C^{(M)}(U) &:= \{f \in C^\infty(U) : \forall \text{ compact } K \subseteq U ~\forall \rho>0: 
\eqref{eq:1} \text{ is bounded} \}, \\
C^{\{M\}}(U) &:= \{f \in C^\infty(U) : \forall \text{ compact } K \subseteq U ~\exists \rho>0: 
\eqref{eq:1} \text{ is bounded} \}.
\end{align*}
The elements of $C^{(M)}(U)$ are said to be of \emph{Beurling type}; those of $C^{\{M\}}(U)$ of \emph{Roumieu type}. 
If $M_k=1$, for all $k$, then $C^{(M)}(U)$ consists of the restrictions to $U$ of the real and imaginary parts of all entire functions, 
while $C^{\{M\}}(U)$ coincides with the ring $C^\omega(U)$ of real analytic functions on $U$. 
We shall also write $C^{[M]}$ and thereby mean that $C^{[M]}$ stands for either $C^{(M)}$ or $C^{\{M\}}$. 

A sequence $M=(M_k)$ is \emph{log-convex} if $k\mapsto\log(M_k)$ is convex, i.e.,
\begin{equation} \label{eq:2} 
M_k^2 \le M_{k-1} \, M_{k+1} \quad \text{ for all } k.
\end{equation}
If $M=(M_k)$ is log-convex, then 
$k\mapsto (M_k/M_0)^{1/k}$ is increasing 
and
\begin{equation} \label{eq:3}
M_l \, M_k\le M_0\,M_{l+k} \quad \text{ for all }l,k\in \mathbb{N}.
\end{equation}
Let us assume $M_0=1$ from now on.
Furthermore, we have that $k\mapsto k!M_k$ is log-convex (since Euler's $\Gamma$-function is so); 
if $M=(M_k)$ satisfies this weaker condition we say that it is \emph{weakly log-convex}.
If $M=(M_k)$ is weakly log-convex, then $C^{[M]}(U)$ is a ring, for all open
subsets $U \subseteq \mathbb{R}^n$. 

If $M=(M_k)$ is log-convex, then (see the proof of \cite[2.9]{KMRc}) we have 
\begin{equation} \label{eq:4}
M_1^j \, M_k\ge M_j\, M_{\alpha_1} \cdots M_{\alpha_j} \quad \text{ for all }
\alpha_i\in \mathbb{N}_{>0}\text{ with } \alpha_1+\dots+\alpha_j = k.
\end{equation}
Condition \eqref{eq:4} implies that the class of $C^{[M]}$-mappings is stable 
under composition. This is due to \cite{Roumieu62/63} in the Roumieu case, see also \cite{Dynkin80} or \cite[4.7]{BM04}; 
the same proof works in the Beurling case. We reproof it in Theorem \ref{nmb:4.11}; compare also with Lemma \ref{nmb:2.5}. For a partial converse, 
see \cite{RainerSchindl12}. 

If $M=(M_k)$ is log-convex, then the inverse function theorem
for $C^{\{M\}}$ holds (\cite{Komatsu79}; see also \cite[4.10]{BM04}), and 
$C^{\{M\}}$ is closed under solving ODEs (due to \cite{Komatsu80}).
If additionally we have $M_{k+1}/M_k \to \infty$, then also $C^{(M)}$ is closed under taking the inverse and solving ODEs
(again by \cite{Komatsu79} and \cite{Komatsu80}). 
See \cite{Yamanaka89}, \cite{Yamanaka91}, \cite{RainerSchindl14}, and Section \ref{nmb:9.2} for Banach space versions of these results.

Suppose that $M=(M_k)$ and $N=(N_k)$ are such that 
$\sup_k(M_k/N_k)^{1/k}<\infty$, i.e.\
\begin{equation} \label{eq:5}
\exists C,\rho>0 ~ \forall k \in \mathbb N : ~M_k \le C \rho^k N_k. 
\end{equation}
Then $C^{(M)}(U) \subseteq C^{(N)}(U)$ and $C^{\{M\}}(U) \subseteq C^{\{N\}}(U)$. 
The converse is true if $M=(M_k)$ is weakly log-convex:
In the Roumieu case the inclusion $C^{\{M\}}(U) \subseteq C^{\{N\}}(U)$ 
implies \eqref{eq:5} thanks to the existence of a function 
$f \in C^{\{M\}}(\mathbb{R})$ such that $|f^{(k)}(0)| \ge k! \, M_k$ for all $k$ (see \cite[Thm.\ 1]{Thilliez08}; and also Section \ref{nmb:2.3}).
In the Beurling case the equivalence of $C^{(M)}(U) \subseteq C^{(N)}(U)$ and \eqref{eq:5} follows from  
the closed graph theorem; see Bruna \cite{Bruna80/81}.
As a consequence we see that the following three conditions are equivalent:
 $C^\omega(U)$ is contained in $C^{\{M\}}(U)$, the restrictions of entire functions are contained in $C^{(M)}(U)$, and $\varliminf M_k^{1/k} > 0$.

$C^{[M]}$ is \emph{stable under derivations} 
(alias \emph{derivation closed}) if 
\begin{equation} \label{eq:6}
\sup_{k \in \mathbb{N}_{>0}} \Big(\frac{M_{k+1}}{M_k}\Big)^{\frac{1}{k}} < \infty.
\end{equation}
The converse is true if $M=(M_k)$ is weakly log-convex: 
$C^{\{M\}}$ is stable under derivations if and only if \eqref{eq:6} holds.

A sequence $M=(M_k)$ is said to have 
\emph{moderate growth} if
\begin{equation} \label{eq:7}
\sup_{j,k \in \mathbb{N}_{>0}} \Big(\frac{M_{j+k}}{M_j \, M_k}\Big)^{\frac{1}{j+k}} <
\infty.
\end{equation}
Moderate growth implies \eqref{eq:6} and thus stability under derivations. 
If $M=(M_k)$ is weakly log-convex and has moderate growth, then $C^{[M]}(U)$ is \emph{stable under ultra\-differential operators},
%An operator of the form $P(D) = \sum_\alpha a_\alpha D^\alpha$, $a_\alpha \in \mathbb{C}$, is an \emph{ultradifferential operator of class $C^{(M)}$ 
%(resp.\ $C^{\{M\}}$)}
%if there are constants $C,\rho>0$ (resp.\ for each $\rho>0$ there exists $C>0$) such that 
%$|a_\alpha| \le \frac{C \rho^{|\alpha|}}{|\alpha|!\,M_{|\alpha|}}$ for all $\alpha$, 
see \cite[2.11 and 2.12]{Komatsu73}.

For sequences $M=(M_k)$ and $N=(N_k)$ of positive real numbers we define
\begin{align*}
M \lhd N \quad :&\Leftrightarrow \quad \forall \rho>0 ~\exists C>0 : M_k \le C \rho^k N_k ~\forall k \in \mathbb{N}\\
&\Leftrightarrow \quad \lim_{k \to \infty} \Big(\frac{M_k}{N_k}\Big)^{\frac{1}{k}}=0. 
\end{align*}
If $M\lhd N$, then we have $C^{\{M\}}(U) \subseteq C^{(N)}(U)$.
If $M=(M_k)$ is weakly log-convex, also the converse is true: $C^{\{M\}}(U) \subseteq C^{(N)}(U)$ implies $M\lhd N$.
This follows from the existence of a function $f \in C^{\{M\}}(\mathbb{R})$ with $|f^{(k)}(0)| \ge k! \, M_k$ for all $k$ 
(see \cite[Thm.\ 1]{Thilliez08}).
As a consequence $C^\omega(U)$ is contained in $C^{(M)}(U)$ if and only if $M_k^{1/k} \to \infty$.

\begin{theorem}[Denjoy--Carleman \cite{Denjoy21}, \cite{Carleman26}]\label{nmb:2.2} 
For a sequence $M=(M_k)$ of positive real numbers 
the following statements are equivalent:
\begin{enumerate}
\item[(1)] $C^{[M]}$ is \emph{quasianalytic}, i.e., for open connected $U
\subseteq \mathbb{R}^n$ and each $x \in U$, the Taylor series homomorphism centered at $x$ 
from $C^{[M]}(U,\mathbb R)$ into the space of 
formal power series is injective.
\item[(2)] $\sum_{k=1}^\infty \frac1{m_k^{\flat(i)}} = \infty$
where $m_k^{\flat(i)}:= \inf\{(j!\,M_j)^{1/j}: j\ge k\}$ is the increasing minorant of 
$(k!\,M_k)^{1/k}$.
\item[(3)] $\sum_{k=1}^\infty (\frac1{M_k^{\flat(lc)}})^{1/k} = \infty$
where $M^{\flat(lc)}_k$ is the log-convex minorant of $k!\,M_k$, given by 
$M^{\flat(lc)}_k := \inf\{(j!\,M_j)^{\frac{l-k}{l-j}}(l!\,M_l)^{\frac{k-j}{l-j}}: j\le k\le l, j<l\}$. 
\item[(4)] 
$\sum_{k=0}^\infty \frac{M^{\flat(lc)}_k}{M^{\flat(lc)}_{k+1}}=\infty$.
\end{enumerate}
\end{theorem}

For contemporary proofs of the equivalence of \thetag{2}, \thetag{3}, \thetag{4} and quasianalyticity of $C^{\{M\}}$,
%\thetag{1} in the Roumieu case $C^{[M]} = C^{\{M\}}$, 
see for instance \cite[1.3.8]{Hoermander83I} 
or \cite[19.11]{Rudin87}. 
For the equivalence of these conditions to the quasianalyticity of $C^{(M)}$, see \cite[4.2]{Komatsu73}.

\subsection{\label{nmb:2.3}Sequence spaces}
Let $M=(M_k)_{k\in\mathbb N}$ be a sequence of positive real numbers, and $\rho>0$.
We consider (where $\mathcal F$ stands for `formal power series')
\begin{gather*}
\mathcal F^{M}_\rho := \Bigl\{(f_k)_{k\in\mathbb N}\in\mathbb R^\mathbb N:
\exists C>0\,\forall k\in\mathbb N:|f_k|\leq C\rho^k k! M_k\Bigr\},\\
\mathcal F^{(M)} := \bigcap_{\rho>0}\mathcal F^M_\rho,\quad \text{ and } \quad 
\mathcal F^{\{M\}} := \bigcup_{\rho>0}\mathcal F^M_\rho.
\end{gather*}

\begin{lemma*} 
Consider the following conditions for two positive sequences $M^i=(M^i_k)$, $i=1,2$, and $0<\sigma<\infty$:
\begin{enumerate}
\item[\thetag{1}] $\sup_k (M^1_k/M^2_k)^{1/k} =\sigma$.
\item[\thetag{2}] For all $\rho>0$ we have $\mathcal{F}^{M^1}_{\rho} \subseteq \mathcal{F}^{M^2}_{\rho \sigma}$.
\item[\thetag{3}] $\mathcal{F}^{\{M^1\}} \subseteq \mathcal{F}^{\{M^2\}}$.
\item[\thetag{4}] $\mathcal{F}^{(M^1)} \subseteq \mathcal{F}^{(M^2)}$.
\item[\thetag{5}] $M^1 \lhd M^2$.
\item[\thetag{6}] $\mathcal{F}^{\{M^1\}} \subseteq \mathcal{F}^{(M^2)}$.
\end{enumerate}
Then we have \thetag{1} $\Leftrightarrow$ \thetag{2} $\Leftrightarrow$ \thetag{3} $\Rightarrow$ \thetag{4} 
and \thetag{5} $\Leftrightarrow$ \thetag{6}.
\end{lemma*}

\begin{demo}{Proof}
\thetag{1} $\Rightarrow$ \thetag{2} 
Let $f=(f_k)\in\mathcal F^{M^1}_{\rho}$, i.e., there is a $C>0$ such that 
$|f_k|\le C \rho^{k}k!M^1_k \le C (\rho \sigma)^{k}k!M^2_k$, for all $k$.
So $f \in \mathcal{F}^{M^2}_{\rho \sigma}$. 

\thetag{2} $\Rightarrow$ \thetag{3} and \thetag{2} $\Rightarrow$ \thetag{4} follow by definition. 

\thetag{3} $\Rightarrow$ \thetag{1}
Let $f_k:=k!M^1_k$. Then
$f=(f_k)\in\mathcal F^{\{M^1\}}\subseteq \mathcal F^{\{M^2\}}$, so there
exists $\rho>0$ such that $k!M^1_k\leq \rho^{k+1}k!M^2_k$ for all $k$.

\thetag{5} $\Rightarrow$ \thetag{6} 
Let $f=(f_k) \in \mathcal{F}^{M^1}_{\rho}$. As $M^1 \lhd M^2$, for each $\sigma>0$ there exists $C>0$ such that 
$|f_k| \le C (\frac{\sigma}{\rho})^k k! M_k$ for all $k$. So $f \in \mathcal{F}^{M^2}_{\frac{\sigma}{\rho}}$ for all $\sigma$.

\thetag{6} $\Rightarrow$ \thetag{5}
Since $(k!M_k) \in \mathcal{F}^{\{M^1\}} \subseteq \mathcal{F}^{(M^2)}$, for each $\rho>0$ there is $C>0$ such that $k!M_k^1 \le C \rho^k k! M_k^2$ 
for all $k$, i.e., 
$M^1 \lhd M^2$. 
\qed\end{demo}

\begin{theorem} \label{nmb:2.4}
Let $M=(M_k)$ be a (weakly) log-convex sequence of positive real numbers. 
Then we have
\begin{equation} \label{eq:8}
\mathcal F^{\{M\}} = \bigcap_{L}\mathcal F^{(L)} = \bigcap_{L}\mathcal F^{\{L\}},
\end{equation}
where the intersections are taken over all (weakly) log-convex $L=(L_k)$ with $M \lhd L$. 
\end{theorem}

\begin{demo}{Proof}
The inclusions $\mathcal F^{\{M\}} \subseteq \bigcap_{L}\mathcal F^{(L)} \subseteq \bigcap_{L}\mathcal F^{\{L\}}$
follow from Lemma \ref{nmb:2.3}. So it remains to prove that $\mathcal F^{\{M\}} \supseteq \bigcap_{L}\mathcal F^{\{L\}}$. 
Let $f=(f_k)\notin \mathcal F^{\{M\}}$, i.e., 
\begin{equation} \label{eq:9}
\varlimsup \Big(\frac{|f_k|}{k! M_k}\Big)^{\frac{1}{k}} = \infty.
\end{equation}
We must show that there exists a (weakly) log-convex $L=(L_k)$ with $M \lhd L$ such that $f\notin \mathcal F^{\{L\}}$.

Choose $a_j, b_j>0$ with $a_j\nearrow \infty$, $b_j\searrow 0$, and 
$a_jb_j \to \infty$.
Now \eqref{eq:9} implies that there exists a strictly increasing sequence $k_j \in \mathbb{N}$ such that
\[
\Big(\frac{|f_{k_j}|}{(k_j)! M_{k_j}}\Big)^{\frac{1}{k_j}} \ge a_j.
\]
Passing to a subsequence we may assume that 
$k_0>0$ and that \[1 < \beta_j:=b_j\, \Big(\frac{|f_{k_j}|}{(k_j)! M_{k_j}}\Big)^{\frac{1}{k_j}}\nearrow \infty.\] 
Passing to a subsequence again we may also get 
\begin{equation} \label{eq:10}
\beta_{j+1} 
\ge \left(\beta_j\right)^{k_j}.
\end{equation}
We define a piecewise affine function $\phi$ by setting  
\begin{align*}
\phi(k) := 
\begin{cases}
0 & \text{ if } k=0,
\\
k_j \log\beta_j & \text{ if } k=k_j,
\\
c_j + d_j k & \text{ for the minimal $j$ with } k \le k_j,
\end{cases}
\end{align*}
where $c_j$ and $d_j$ are chosen such that $\phi$ is well defined 
and $\phi(k_{j-1})=c_j+d_j k_{j-1}$, i.e., for $j\ge 1$,
\begin{align*}
c_j+d_j k_j &= k_j \log\beta_j,\\
c_j+d_j k_{j-1} &= k_{j-1} \log \beta_{j-1},
\quad \text{ and } \\
c_0 &=0,\\
d_0 &= \log\beta_0. 
\end{align*}
This implies first that $c_j\le 0$ and then
\begin{align*} 
\log\beta_j \le d_j&=\frac{k_j\log \beta_j-k_{j-1}\log\beta_{j-1}}{k_j-k_{j-1}} 
\le \frac{k_j}{k_j-k_{j-1}} \log\beta_j
\\&
\overset{\eqref{eq:10}}\leq\frac{\log\beta_{j+1}}{k_j-k_{j-1}}\leq\log\beta_{j+1}.
\end{align*}
Thus $j\mapsto d_j$ is increasing and so $\phi$ is convex. 
The fact that all $c_j \le 0$ implies that $\phi(k)/k$ is increasing.

Now let
\[
L_k := e^{\phi(k)} \cdot M_k.
\]
Then $L=(L_k)$ is (weakly) log-convex, since so is $M=(M_k)$. %
As $\phi(k)/k$ is increasing and 
$e^{\phi(k_j)/k_j} = \beta_j \to \infty$, we find $M \lhd L$.
Finally,
$f\notin \mathcal F^{\{L\}}$, since we have 
\[
\Big(\frac{|f_{k_j}|}{(k_j)! L_{k_j}}\Big)^{\frac{1}{k_j}} = \Big(\frac{|f_{k_j}|}{(k_j)! M_{k_j}}\Big)^{\frac{1}{k_j}} \cdot e^{-\phi(k_j)/k_j} = \Big(\frac{|f_{k_j}|}{(k_j)! M_{k_j}}\Big)^{\frac{1}{k_j}} \cdot \beta_j^{-1} = b_j^{-1} \to \infty.
\]
The proof is complete.
\qed\end{demo}

\begin{remark*}
\thetag{1} If $M_0 =1 \le M_1$ we also have $L_0 =1 \le L_1$.

\thetag{2} The proof also shows that, if $M=(M_k)$ is just any positive sequence, then 
\eqref{eq:8} still holds if the 
intersections are taken over all positive sequences $L=(L_k)$ with $M \lhd L$. 
\end{remark*}

\begin{lemma}\label{nmb:2.5}
Let $M=(M_k)$ and $L=(L_k)$ be sequences of positive real numbers.
Then for the composition of formal power series we have
\begin{equation} \label{eq:11}
\mathcal F^{[M]}\o \mathcal F^{[L]}_{>0}\subseteq \mathcal F^{[M\o L]},
\end{equation}
where $(M\o L)_k := \max\{M_jL_{\alpha_1}\dots L_{\alpha_j}: \alpha_i\in \mathbb{N}_{>0}, \alpha_1+\dots+\alpha_j = k \}$. 
\end{lemma}

Here $\mathcal F^{[L]}_{>0}$ 
is the space of formal power series in $\mathcal F^{[L]}$ with vanishing constant term.

\begin{demo}{Proof} 
Let $f \in \mathcal F^{(M)}$ and $g \in \mathcal F^{(L)}$ (resp.\ $f \in \mathcal F^{\{M\}}$ and $g \in \mathcal F^{\{L\}}$).
For $k>0$ we have (inspired by \cite{FaadiBruno1855})
\begin{align*}
\frac{(f\o g)_k}{k!} :&= \sum_{j=1}^k \frac{f_j}{j!}
\sum_{\substack{\alpha\in \mathbb{N}_{>0}^j\\ \alpha_1+\dots+\alpha_j =k}}
\frac{g_{\alpha_1}}{\alpha_1!}\dots
\frac{g_{\alpha_j}}{\alpha_j!},
\\
\frac{|(f\o g)_k|}{k!(M\o L)_k} 
&\le \sum_{j=1}^k \frac{|f_j|}{j!M_j}
\sum_{\substack{\alpha\in \mathbb{N}_{>0}^j\\ \alpha_1+\dots+\alpha_j =k}}
\frac{|g_{\alpha_1}|}{\alpha_1!L_{\alpha_1}}\dots
\frac{|g_{\alpha_j}|}{\alpha_j!L_{\alpha_j}}
\\&
\le \sum_{j=1}^k \rho_f^{j}C_f
\sum_{\substack{\alpha\in \mathbb{N}_{>0}^j\\ \alpha_1+\dots+\alpha_j =k}}
\rho_g^k C_g^j
\le \sum_{j=1}^k \rho_f^{j}C_f
\binom{k-1}{j-1}
\rho_g^k C_g^j
\\&
= \rho_g^k\rho_f C_f C_g \sum_{j=1}^k (\rho_f C_g)^{j-1}
\binom{k-1}{j-1}
= \rho_g^k\rho_f C_f C_g (1+\rho_f C_g)^{k-1}
\\&
= (\rho_g(1+\rho_f C_g))^k\frac{\rho_f C_f C_g}{1+\rho_f C_g}.
\end{align*}
This implies \eqref{eq:11} in the Roumieu case.
For the Beurling case,
let $\tau > 0$ be arbitrary, and choose $\sigma>0$ such that $\tau=\sqrt{\sigma}+\sigma$. 
If we set $\rho_g = \sqrt{\sigma}$ and $\rho_f = \sqrt{\sigma}/C_g$, then $f \o g \in \mathcal{F}^{M \o L}_{\tau}$.
\qed\end{demo}

\subsection{Convention} \label{nmb:2.6}
For a positive sequence $M = (M_k) \in (\mathbb{R}_{>0})^\mathbb{N}$
consider the following properties:
\begin{enumerate}
\item[\thetag{0}] $M_0=1\le M_1$.
\item $M=(M_k)$ is weakly log-convex, i.e., $k \mapsto \log(k!\, M_k)$ is convex.
\item $M=(M_k)$ is log-convex, i.e., $k \mapsto \log(M_k)$ is convex.
\item $M=(M_k)$ is derivation closed, i.e., $k \mapsto (\frac{M_{k+1}}{M_k})^{\frac{1}{k}}$ is bounded.
\item $M=(M_k)$ has moderate growth, i.e., $(j,k) \mapsto (\frac{M_{j+k}}{M_j\, M_k})^{\frac{1}{j+k}}$ is bounded.
\item $\frac{M_{k+1}}{M_k} \to \infty$. 
\item $M_k^{1/k} \to \infty$, or equivalently, $C^\omega \subseteq C^{(M)}$.
\end{enumerate}
Henceforth, if not specified otherwise, we assume that $M=(M_k)$, $N=(N_k)$, $L=(L_k)$, etc., satisfy condition \thetag{0}.
It will be explicitly stated when some of the other 
properties \thetag{1}--\thetag{6} are assumed.

\begin{remarks*}
Let $M=(M_k)$ be a positive sequence. 
We may replace $(M_k)_k$ by $(C \rho^k M_k)_k$ with $C,\rho>0$
without changing $\mathcal F^{[M]}$ (see Section \ref{nmb:2.3}). In particular, it is no loss of generality
to assume that $M_1>1$ (put $C\rho>1/M_1$) and $M_0=1$ (put $C:=1/M_0$).
Each one of the properties \thetag{1}--\thetag{6} is preserved by this modification.
Furthermore $M=(M_k)$ is quasianalytic if and only if the modified sequence is so, since $(M_k^{\flat(lc)})_k$ (see Theorem \ref{nmb:2.2}) 
is modified in the same way.

Conditions \thetag{0} and \thetag{1} together imply that $k\mapsto k!\, M_k$ is monotone increasing,
while \thetag{0} and \thetag{2} together imply that $k\mapsto M_k$ is monotone increasing.
\end{remarks*}

\section{\label{nmb:3}Whitney jets on Banach spaces}

\subsection{\label{nmb:3.1}Whitney jets}
Let $E$ and $F$ be Banach spaces.
For open $U\subseteq E$ consider the space $C^\infty(U,F)$ of 
arbitrarily often Fr\'echet differentiable mappings $f:U\to F$.
For such $f$ we have the derivatives $f^{(k)}:U\to L_{\text{sym}}^k(E,F)$,
where $L_{\text{sym}}^k(E,F)$ denotes the space of symmetric $k$-linear bounded mappings $E\times\dots\times E\to F$.
We also have the iterated uni-directional derivatives $d^k_vf(x)\in F$ defined by 
\[
d^k_v f(x):=\left(\frac{d}{dt}\right)^k f(x+t\,v)|_{t=0}.
\]
Let $j^\infty:C^\infty(U,F)\to J^\infty(U,F):=\prod_{k\in\mathbb N}C(U,L_{\text{sym}}^k(E,F))$ be the jet mapping
$f\mapsto (f^{(k)})_{k\in\mathbb N}$.
On $L_{\text{sym}}^k(E,F)$ we consider the operator norm
\[
\|\ell\|_{L_{\text{sym}}^k(E,F)} :=\sup\Bigl\{\|\ell(v_1,\dots,v_k)\|_F:\|v_j\|_E\leq 1\text{ for all 
}j\in\{1,\dots,k\}\Bigr\}.
\]
Note that by the polarization equality (see \cite[7.13.1]{KM97}) 
\begin{align*}
\sup\{\|\ell(v,\dots,v)\|_F:\|v\|_E\leq 1\} 
&\leq \|\ell\|_{L_{\text{sym}}^k(E,F)}
\\&
\leq (2e)^k\,\sup\{\|\ell(v,\dots,v)\|_F:\|v\|_E\leq 1\} 
\end{align*}
For an infinite jet $f=(f^k)_{k\in\mathbb N}\in \prod_{k\in\mathbb N}L_{\text{sym}}^k(E,F)^X$ on a subset $X\subseteq E$
let the Taylor polynomial $(T_y^n f)^k:X\to L_{\text{sym}}^k(E,F)$ of order $n$ at $y$ be
\[
(T_y^n f)^k(x)(v_1,\dots,v_k) :=\sum_{j=0}^n \frac1{j!}\,f^{j+k}(y)(x-y,\dots,x-y,v_1,\dots v_k)
\]
and the remainder
\[
(R_y^n f)^k(x) := f^k(x) - (T_y^n f)^k(x)=(T_x^n f)^k(x)-(T_y^n f)^k(x) \in L_{\text{sym}}^k(E,F).
\]

Let
\begin{align*}
\|f\|_{k} &:= \sup\{\|f^k(x)\|_{L_{\text{sym}}^k(E,F)}:x\in X\}\in[0,+\infty]\quad \text{and}\\
|||f|||_{n,k} 
&:=\sup\Bigl\{(n+1)!\frac{\|(R_y^nf)^k(x)\|_{L_{\text{sym}}^k(E,F)}}{\|x-y\|^{n+1}}:x,y\in X,x\ne y\Bigr\}\in[0,+\infty].
\end{align*}
By Taylor's theorem, 
for $f\in C^\infty(U,F)$ and 
$[x,y]\subseteq U$ we have 
\begin{align*}
(R_y^nf)^k(x)
&=f^{(k)}(x) - \sum_{j \leq n} \frac{f^{(k+j)}(y)(x-y)^j}{j!} \\
&= \int_0^1 \frac{(1-t)^n}{n!}\, f^{(k+n+1)}(y+t(x-y))(x-y)^{n+1}\,dt
\end{align*}
and hence for convex $X \subseteq U$:
\begin{align}
|||&j^\infty f|_X|||_{n,k} := \notag
\\&
=\sup\Bigl\{(n+1)!\frac{\|(R_y^nf)^k(x)(v_1,\dots,v_k)\|_F}{\|x-y\|^{n+1}}:
\|v_j\|_E\leq 1,x,y\in X,x\ne y\Bigr\} \notag
\\ 
&\leq\sup\Bigl\{\frac{\|f^{(k+n+1)}(x)(v_1,\dots,v_k,x-y,\dots,x-y)\|_F}{\|x-y\|^{n+1}}:
\|v_j\|_E\leq 1,x\ne y\Bigr\} \notag
\\ 
&\leq \|j^\infty f|_X\|_{n+k+1}.
\label{eq:12}
\end{align}
We supply $C^\infty(U,F)$ with the semi-norms
\[
f\mapsto \|j^\infty f|_K\|_n\quad \text{ for all compact }K\subseteq U\text{ and all }n\in\mathbb N.
\]
For compact convex $K\subseteq E$ the space $C^\infty(E\supseteq K,F)$ of Whitney jets on $K$
is defined by 
\begin{multline*}
C^\infty(E\supseteq K,F):=
\\
=\Bigl\{f=(f^k)_{k\in\mathbb N}\in\prod_{k\in\mathbb N} C(K,L_{\text{sym}}^k(E,F)):|||f|||_{n,k}<\infty
\text{ for all }n,k\in\mathbb N\Bigr\}
\end{multline*}
and is supplied with the seminorms $\|~\|_n$ for $n\in\mathbb N$
together with $|||~|||_{n,k}$ for $n,k\in\mathbb N$.

\begin{lemma} \label{nmb:3.2}
	For Banach spaces $E$ and $F$ and compact convex $K\subseteq E$ the space $C^\infty(E\supseteq K,F)$
	is a Fr\'echet space.	
\end{lemma}

\begin{demo}{Proof}
The injection of $C^\infty(E\supseteq K,F)$ into $\prod_{k\in\mathbb N}C(K,L_{\text{sym}}^k(E,F))$
is continuous by definition and $C(K,L_{\text{sym}}^k(E,F))$ is a Banach space, 
so a Cauchy sequence $(f_p)_p$ in $C^\infty(E\supseteq K,F)$ has an infinite jet 
$f_\infty=(f_\infty^k)_k$ as component-wise limit in $\prod_{k\in\mathbb N}C(K,L_{\text{sym}}^k(E,F))$ with respect to the 
seminorms $\|~\|_n$.
This is the limit also with respect to the finer structure of $C^\infty(E\supseteq K,F)$ with the 
additional seminorms $|||~|||_{n,k}$ as follows: 
For given $n,k\in\mathbb N$ and $\epsilon>0$ there exists by the Cauchy condition a $p_0$ such that 
$|||f_p-f_q|||_{n,k}<\epsilon/3$ for all $p,q\geq p_0$.
By the convergence $f_q\to f_\infty$ in $\prod_{k\in\mathbb N}C(K,L_{\text{sym}}^k(E,F))$ there exists for given $x,y\in K$
with $x\ne y$
a $q\geq p_0$ such that for all $m\leq k+n$ 
\[
\|f_q-f_\infty\|_m\leq \frac{\|x-y\|^{n+1}}{(n+1)!}\,\frac{\epsilon}3\, \min\{1,e^{-\|x-y\|}\}
\]
and hence 
\begin{align*}
\|(T^n_y f_q)^k(x) &- (T^n_y f_\infty)^k(x)\|_{L_{\text{sym}}^k(E,F)} \leq
\\&
\leq \sum_{j=0}^n \|f_q^{k+j}(y)-f_\infty^{k+j}(y)\|_{L_{\text{sym}}^{k+j}(E,F)}\,\frac{\|x-y\|^j}{j!} \\
&\leq \sum_{j=0}^n \|f_q-f_\infty\|_{k+j}\,\frac{\|x-y\|^j}{j!} \\
&\leq \frac{\|x-y\|^{n+1}}{(n+1)!}\,\frac{\epsilon}{3}\,e^{-\|x-y\|}\, \sum_{j=0}^n \frac{\|x-y\|^j}{j!} 
\leq \frac{\|x-y\|^{n+1}}{(n+1)!}\,\frac{\epsilon}{3}.
\end{align*}
So 
\begin{align*}
(n+1)!\, &\frac{\|(R^n_y f_p)^k(x)-(R^n_y f_\infty)^k(x)\|_{L_{\text{sym}}^k(E,F)}}{\|x-y\|^{n+1}}\leq \\
&\leq |||f_p-f_q|||_{n,k} 
+(n+1)!\, \frac{\|(R^n_y f_q)^k(x)-(R^n_y f_\infty)^k(x)\|_{L_{\text{sym}}^k(E,F)}}{\|x-y\|^{n+1}} \\
&\leq \frac{\epsilon}3
+(n+1)!\, \frac{\|f_q^k(x)-f_\infty^k(x)\|_{L_{\text{sym}}^k(E,F)}}{\|x-y\|^{n+1}} \\
&\qquad +(n+1)!\, \frac{\|(T^n_y f_q)^k(x)-(T^n_y f_\infty)^k(x)\|_{L_{\text{sym}}^k(E,F)}}{\|x-y\|^{n+1}} \\
&\leq 3\,\frac{\epsilon}3=\epsilon
\end{align*}
and finally
\begin{align*}
|||f_p-f_\infty|||_{n,k} \leq \epsilon\text{ for all }p\geq p_0.
\end{align*}
Consequently,
\[
|||f_\infty|||_{n,k}\leq |||f_\infty-f_p|||_{n,k} + |||f_p|||_{n,k}<\infty,
\]
i.e.,
$f_\infty\in C^\infty(E\supseteq K,F)$.
\qed\end{demo}

\section{\label{nmb:4}The category of Denjoy--Carleman differentiable mappings}

\subsection{\label{nmb:4.1}Spaces of Denjoy--Carleman jets and mappings between Banach spaces}
Let %
$E$ and $F$ be Banach spaces, $K\subseteq E$ compact, and $\rho>0$.
Let
\begin{align*}
C^{M}_\rho(E\supseteq K,F) :&= 
\Bigl\{(f^m)_m\in \prod_{m\in\mathbb N}C(K,L_{\text{sym}}^m(E,F)):\|f\|_\rho<\infty\Bigr\},\\
\intertext{where }\|f\|_\rho:&=\max\biggl\{\sup\Bigl\{\frac{\|f\|_m}{m!\rho^m M_m}:m\in\mathbb N\Bigr\},\\
&\qquad\qquad\sup\Bigl\{\frac{|||f|||_{n,k}}{(n+k+1)!\,\rho^{n+k+1}\,M_{n+k+1}}:k,n\in\mathbb N\Bigr\}\biggr\},\\
\end{align*}
cf.\  \cite[11]{ChaumatChollet94}, \cite[11]{ChaumatChollet98} 
and \cite[3]{SchmetsValdivia03}, and,
for an open neighborhood $U$ of $K$ in $E$, let
\begin{align*}
C^M_{K,\rho}(U,F):&=\Bigl\{f\in C^\infty(U,F):j^\infty f|_K\in C^M_\rho(E\supseteq K,F)\Bigr\}
\end{align*}
supplied with the semi-norm $f\mapsto \|j^\infty f|_K\|_\rho$.
This space is not Hausdorff and for infinite dimensional $E$ it(s Hausdorff quotient) will not always
be complete. This is the reason for considering the jet spaces $C^{M}_\rho(E\supseteq K,F)$ instead.
Note that for convex $K$ we have $|||j^\infty f|_K|||_{n,k}\leq \|j^\infty f|_K\|_{n+k+1}$ by 
\eqref{eq:12} 
and hence the seminorm $f\mapsto\|j^\infty f|_K\|_\rho$ 
on $C^M_{K,\rho}(U,F)$ coincides with
\[
f\mapsto \sup\Bigl\{\frac{\|f^{(n)}(x)\|_{L_{\text{sym}}^n(E,F)}}{n!\rho^n M_n}:x\in K,n\in\mathbb N\Bigr\} =: \|f\|_{K,\rho}.
\]
Thus
\[
C^M_{K,\rho}(U,F)=\Bigl\{f\in C^\infty(U,F):(\|j^\infty f|_K\|_m)_m\in \mathcal F^M_\rho\Bigr\}
\]
and the bounded subsets $\mathcal{B}\subseteq C^M_{K,\rho}(U,F)$ are exactly those $\mathcal{B}\subseteq C^\infty(U,F)$
for which $(b_m)_m\in \mathcal F^M_\rho$, where $b_m:=\sup\{\|j^\infty f|_K\|_m:f\in\mathcal{B}\}$.

For open convex $U\subseteq E$ and compact convex $K \subseteq U$ let
\begin{align*}
C^{(M)}(E\supseteq K,F):&=\bigcap_{\rho>0} C^{M}_\rho(E\supseteq K,F), \\
C^{\{M\}}(E\supseteq K,F):&=\bigcup_{\rho>0} C^{M}_\rho(E\supseteq K,F),\text{ and} \\
C^{[M]}(U,F):&=\Bigl\{f\in C^\infty(U,F):\forall K: (f^{(k)}|_K)\in C^{[M]}(E\supseteq K,F)\Bigr\}.
\end{align*}
That means, we consider
the projective limit
\begin{align*}
C^{(M)}(E\supseteq K,F)&:=\varprojlim_{\rho>0} C^M_\rho(E\supseteq K,F),\\
\end{align*}
the inductive limit
\begin{align*}
C^{\{M\}}(E\supseteq K,F)&:=\varinjlim_{\rho>0} C^M_\rho(E\supseteq K,F),\\
\end{align*}
and the projective limits
\begin{align*}
C^{[M]}(U,F)&:=\varprojlim_{K\subseteq U} C^{[M]}(E\supseteq K,F),
\end{align*}
where $K$ runs through all compact convex subsets of $U$.

Furthermore, we consider
the projective limit
\begin{align*}
C^{(M)}_K(U,F)&:=\varprojlim_{\rho>0} C^M_{K,\rho}(U,F),
\end{align*}
and the inductive limit
\begin{align*}
C^{\{M\}}_K(U,F)&:=\varinjlim_{\rho>0} C^M_{K,\rho}(U,F).
\end{align*}
Thus

\[
C^{[M]}_K(U,F)=\Bigl\{f\in C^\infty(U,F):(\|j^\infty f|_K\|_m)_m\in\mathcal F^{[M]}\Bigr\}.
\]
Furthermore, the bounded subsets $\mathcal{B}\subseteq C^{[M]}_K(U,F)$ are exactly those $\mathcal{B} \subseteq C^\infty(U,F)$
for which $(b_m)_m\in \mathcal F^{[M]}$, where $b_m:=\sup\{\|j^\infty f|_K\|_m:f\in\mathcal{B}\}$.
%and subsets 
%$\mathcal{B}\subseteq C^\infty(U,F)$
%for which $(b_m)_m\in \mathcal F^{\{M\}}$, where $b_m:=\sup\{\|j^\infty f|_K\|_m:f\in\mathcal{B}\}$,
%are bounded in $C^{\{M\}}_K(U,F)$.

Finally, the projective limits
\begin{align*}
\varprojlim_{K\subseteq U} C^{[M]}_K(U,F)=\Bigl\{f\in C^\infty(U,F):\forall K:(\|j^\infty f|_K\|_m)_m\in\mathcal F^{[M]}\Bigr\},
\end{align*}
where $K$ runs through all compact convex subsets of $U$, are for 
$E=\mathbb R^n$ and 
$F=\mathbb R$ exactly the vector spaces of Section \ref{nmb:2.1} and the topology is the usual one.

For the inductive limits with respect to $\rho>0$ it suffices to take $\rho \in \mathbb N$ only.

\begin{proposition} \label{nmb:4.2}
	We have the following completeness properties:
\begin{itemize}
\item[(1)] The spaces $C^{M}_\rho(E\supseteq K,F)$ are Banach spaces.
\item[(2)] The spaces $C^{(M)}(E\supseteq K,F)$ are Fr\'echet spaces.
\item[(3)] The spaces $C^{\{M\}}(E\supseteq K,F)$ are compactly regular (i.e., 
compact subsets are contained and compact in some step) (LB)-spaces
hence {($c^\infty$-)}complete, webbed and (ultra-)bornological.
\item[(4)] The spaces $C^{[M]}(U,F)$ are complete spaces.
\item[(5)] As locally convex spaces 
\[
C^{[M]}(U,F):=\varprojlim_{K\subseteq U} C^{[M]}(E\supseteq K,F)=\varprojlim_{K\subseteq U} C^{[M]}_K(U,F).
\]
\end{itemize}	
\end{proposition}

\begin{demo}{Proof}
\thetag{1}
The injection $C^{M}_\rho(E\supseteq K,F)\to \prod_{k\in\mathbb N}C(K,L_{\text{sym}}^k(E,F))$
is by definition continuous and $C(K,L_{\text{sym}}^k(E,F))$ is a Banach space, 
so a Cauchy sequence $(f_p)_p$ in $C^{M}_\rho(E\supseteq K,F)$ has an infinite jet 
$f_\infty=(f_\infty^k)_k$ as component-wise limit in $\prod_{k\in\mathbb N}C(K,L_{\text{sym}}^k(E,F))$.
This is the limit also with respect to the finer structure of $C^{M}_\rho(E\supseteq K,F)$ as follows: 
For fixed $n,k$ and $x\ne y$ we have that
$(R_y^n f_p)^k(x)$ converges to $(R_y^n f_\infty)^k(x)$. So 
we choose for $\epsilon>0$ a $p_0\in\mathbb N$ such that $\|f_p-f_q\|_\rho<\epsilon/2$ for all $p,q\geq p_0$
and
given $x,y,n$, and $k$
we can choose $q>p_0$ such that 
\begin{align*}
(n+1)!\,\frac{\|(R_y^nf_q)^k(x)-(R_y^nf_\infty)^k(x)\|_{L_{\text{sym}}^k(E,F)}}{(n+k+1)!\rho^{n+k+1}M_{n+k+1}\|x-y\|^{n+1}}
&<\frac{\epsilon}2 \\
\intertext{and}
\frac{\|f_q^n(x)-f_\infty^n(x)\|_{L_{\text{sym}}^n(E,F)}}{n!\,\rho^n\,M_n}
&< \frac{\epsilon}2.
\end{align*}
Thus 
\begin{align*}
(n+1)!\,&\frac{\|(R_y^nf_p)^k(x)-(R_y^nf_\infty)^k(x)\|_{L_{\text{sym}}^k(E,F)}}{(n+k+1)!\rho^{n+k+1}M_{n+k+1}\|x-y\|^{n+1}}<\\
&\qquad\qquad\|f_p-f_q\|_\rho + 
(n+1)!\,\frac{\|(R_y^nf_q)^k(x)-(R_y^nf_\infty)^k(x)\|_{L_{\text{sym}}^k(E,F)}}{(n+k+1)!\rho^{n+k+1}M_{n+k+1}\|x-y\|^{n+1}}
<\epsilon
\end{align*}
and hence 
\[
\frac{|||f_p-f_\infty|||_{n,k}}{(n+k+1)!\rho^{n+k+1}M_{n+k+1}}\leq \epsilon
\]
and similarly for
$\frac{\|f_p-f_\infty\|_n}{n!\,\rho^n\,M_n}$.
Thus $\|f_p-f_\infty\|_\rho\leq\epsilon$ for all $p\geq p_0$.

\thetag{2} This is obvious; they are countable projective limits of Banach spaces.

\thetag{3} 
For finite dimensional $E$ and $F$ it is shown in 
\cite{SchmetsValdivia03} that the connecting mappings are nuclear.
For infinite dimensional $E$
the connecting mappings in $C^{\{M\}}(E\supseteq K,F)=\varinjlim_{\rho>0}C^M_\rho(E\supseteq K,F)$
cannot be compact, since the set $\{\ell\in E':\|\ell\|\leq 1\}$ is 
bounded in $C^M_\rho(E\supseteq K,\mathbb R)$ for each $\rho\geq 1$.
In fact, $\|\ell\|_0=\sup\{|\ell(x)|:x\in K\}\leq \sup\{\|x\|:x\in K\}$,
$\|\ell\|_1=\|\ell\|\leq 1$ and $\|\ell\|_m=0$ for $m\geq 2$,
moreover, $(R^n_y\ell)^k=0$ for $n+k\geq 1$ and $(R^0_y\ell)^0(x)=\ell(x-y)$.
It is not relatively compact in any of the spaces $C^M_\rho(E\supseteq K,\mathbb R)$, $\rho\geq 1$, since it is not even pointwise relatively compact
in $C(K,L(E,\mathbb R))$.

In order to show that the (LB)-space in \thetag{3} is compactly regular 
it suffices by \cite[Satz 1]{Neus78} to verify condition (M) of \cite{Retakh70}:
There exists a sequence of increasing 0-neighborhoods $U_n\subseteq C^M_n(E\supseteq K,F)$, such that for 
each $n$ there exists an $m\geq n$ for which the topologies of $C^M_k(E\supseteq K,F)$ and
of $C^M_m(E\supseteq K,F)$ coincide on $U_n$ for all $k\geq m$.

For $\rho'\geq\rho$ we have $\|f\|_{\rho'}\leq\|f\|_\rho$. 
So consider the $\epsilon$-balls $U^\rho_\epsilon(f):=\{g:\|g-f\|_\rho\leq \epsilon\}$ in $C^M_\rho(E\supseteq K,F)$.
It suffices to show that for $\rho>0$, $\rho_1:=2\rho$, $\rho_2>\rho_1$, $\epsilon>0$, and 
$f\in U^\rho_1:=U^\rho_1(0)$ there exists a $\delta>0$ such that 
$U^{\rho_2}_\delta(f)\cap U^{\rho}_1\subseteq U^{\rho_1}_\epsilon(f)$.
Since $f\in U^\rho_1$ we have
\[
\|f\|_n\leq n!\rho^n M_n\text{ and }
|||f|||_{n,k}\leq (n+k+1)!\rho^{n+k+1} M_{n+k+1}\text{ for all }n,k.
\]
Let $\frac1{2^N}<\frac{\epsilon}2$ and $\delta:=\epsilon\,\left(\frac{\rho_1}{\rho_2}\right)^{N-1}$.
Let $g\in U^{\rho_2}_\delta(f)\cap U^{\rho}_1$, i.e.,
\begin{align*}
\|g\|_n&\leq n!\rho^n M_n\text{ for all }n,\\
\|g-f\|_n&\leq \delta\,n!\rho_2^n M_n\text{ for all }n,\\
|||g|||_{n,k}&\leq (n+k+1)!\rho^{n+k+1} M_{n+k+1}\text{ for all }n,k,\\
|||g-f|||_{n,k}&\leq \delta\,(n+k+1)!\rho_2^{n+k+1} M_{n+k+1}\text{ for all }n,k.\\
\end{align*}
Then
\begin{align*}
\|g-f\|_n&\leq \|g\|_n+\|f\|_n\leq 2\,n!\rho^n M_n
=2\,n!\rho_1^n M_n\,\frac1{2^n} \\
&< \epsilon\,n!\rho_1^n M_n\quad\text{ for }n\geq N\\
\intertext{and}
\|g-f\|_n&\leq \delta\,n!\rho_2^n M_n
\leq \epsilon\,n!\rho_1^n M_n\quad\text{ for }n<N. \\
\intertext{Moreover,}
|||g-f|||_{n,k}&\leq 
|||g|||_{n,k}+|||f|||_{n,k}\leq 
2\,(n+k+1)!\rho^{n+k+1} M_{n+k+1} \\
&=2\,(n+k+1)!\rho_1^{n+k+1} M_{n+k+1}\,\frac1{2^{n+k+1}}\\
&< \epsilon\,(n+k+1)!\rho_1^{n+k+1} M_{n+k+1}
\quad\text{ for }n+k+1\geq N \\
\intertext{and}
|||g-f|||_{n,k}&\leq \delta\,(n+k+1)!\rho_2^{n+k+1} M_{n+k+1} \\
&\leq \epsilon\,(n+k+1)!\rho_1^{n+k+1} M_{n+k+1}\quad\text{ for }n+k+1<N.
\end{align*}

\thetag{4} This is obvious; they are projective limits of complete spaces.

\thetag{5} 
Since $j^\infty|_K:C^M_{K,\rho}(U,F)\to C^M_\rho(E\supseteq K,F)$ is by definition a well-defined continuous
linear mapping, it induces such mappings $C^{[M]}_K(U,F)\to C^{[M]}(E\supseteq K,F)$
and $\varprojlim_K C^{[M]}_K(U,F)\to \varprojlim_K C^{[M]}(E\supseteq K,F)$.
The last mapping is obviously injective (use $K:=\{x\}$ for the points $x\in U$).

Conversely, let $f^k_K\in C(K,L_{\text{sym}}^k(E,F))$ be given, such that for each $K$ there exists $\rho>0$ (resp.\ each 
$\rho>0$) we have
$(f^k_K)_{k\in\mathbb N}\in C^M_\rho(E\supseteq K,F)$
and such that $f^k_K|_{K'}=f^k_{K'}$. They define an infinite jet $(f^k)_{k\in\mathbb N}\in J^\infty(U,F)$
by setting $f^k(x):=f^k_{\{x\}}(x)$ which satisfies 
$f^k|_K=f^k_K$ for all $k\in\mathbb N$ and all $K$.

We claim that $f^0\in C^\infty(U,F)$ and $(f^0)^{(k)}=f^k$ for all $k$, i.e.,
$j^\infty f^0|_K=(f^k_K)_k$ for all $k\in\mathbb N$ and all $K$.
\newline
By \cite[5.20]{KM97}
it is enough to show by induction that 
$d^k_vf^0(x)=f^k(x)(v,\dots,v)$.
For $k=0$ this is obvious, so let $k>0$.
Then
\begin{align*}
d^k_v f^0(x)
:&=\lim_{t\to 0}\frac{d^{k-1}_v f^0(x+t\,v)-d^{k-1}_v f^0(x)}{t} \\
&=\lim_{t\to 0}\frac{f^{k-1}(x+t\,v)(v^{k-1})-f^{k-1}(x)(v^{k-1})}{t} \\
&=\lim_{t\to 0}\frac{(R_x^1 f)^{k-1}(x+t\,v)(v^{k-1})}{t} + f^k(x)(v^k) =f^k(x)(v^k).
\end{align*}
Finally, $f^0$ defines an element in $\varprojlim_K C^{[M]}_K(U,F)$, since $\forall K$ we have 
$f^0\in C^M_{K,\rho}(U,F)=\{g\in C^\infty(U,F):j^\infty g|_K\in C^M_\rho(E\supseteq K,F)\}$ for some (resp.\ all) $\rho>0$.

That this bijection is an isomorphism follows, since the seminorm $\|~\|_{K,\rho}$ on $C^M_{K,\rho}(U,F)$ is the 
pull-back of the seminorm $\|~\|_\rho$ on $C^M_\rho(E\supseteq K,F)$. 
\qed\end{demo}

\begin{comment}

\begin{remark}{ nmb.{4.3} Remark}
\begin{itemize}
\item 
We claim regularity of the limit $C^{\{M\}}_K(U,F) =\varinjlim_{\rho} C^M_{K,\rho}(U,F)$
in \cite[1.10]{KMRd}.

\item 
We need regularity in \cite[3.3]{KMRd} for:
\newline
Thus, the mapping $f^\vee : U_1 \to \varinjlim_{\rho_2} C^{M}_{K_2,\rho_2}(U_2,\mathbb{R})$ is $C_b^{\{M\}}$, 
by \cite[2.1R]{KMRq}.
The inductive limit is countable and regular, by \cite[7.4 and 7.5]{Floret71} or \cite[52.37]{KM97}, since it is a Silva space, see
\cite{Komatsu73} (with the same proof in Banach spaces).
Hence, for each compact $K_1 \subseteq U_1$ there exists $\rho_1 >0$ so that the set \thetag{nmb|{expeq2}} is contained and bounded in
$C^{M}_{K_2,\rho_2}(U_2,\mathbb{R})$ (or, equivalently, \thetag{nmb|{expeq3}}), for some $\rho_2>0$.

\item 
And we claim in \cite[4.6R]{KMRd} (cf. \cite[4.6]{KMRq}) that also the
injective limits in
\begin{align*}
\varprojlim_K \varinjlim_{\rho,W} C^M_\rho(W,F) \tag{1} \\
\varprojlim_K \varinjlim_\rho C^M_\rho(K,F) \tag{2} \\
\end{align*}
are regular. 

\item \cite[Theorem 4a]{Bierstedt} weakly compact injective inductive sequences are regular. 
\item \cite[p77]{Bierstedt} An (LF)-space is Mackey-complete if and only if it is regular.
\end{itemize}
\end{remark}
\end{comment}

\subsection{\label{nmb:4.3}Spaces of Denjoy--Carleman differentiable mappings between convenient vector spaces}
For convenient vector spaces $E$ and $F$, and $c^\infty$-open $U\subseteq E$, we define:
\begin{align*}
C_b^{(M)}(U,F) &:= 
\Bigl\{f\in C^\infty(U,F):
\forall B\; \forall\text{ compact } K\subseteq U_B\;\forall \rho>0:\\
&\qquad\big\{\frac{f^{(k)}(x)(v_1,\dots,v_k)}{k!\,\rho^k\, M_k}:k\in \mathbb N,x\in K,\|v_i\|_B\leq 1\big\}\text{ 
is bounded in $F$}\Bigr\}
\\&
= \Bigl\{f\in C^\infty(U,F):
\forall B\; \forall\text{ compact } K\subseteq U_B\;\forall \rho>0:\\
&\qquad\big\{\frac{d^k_vf(x)}{k!\,\rho^k\, M_k}:k\in \mathbb N,x\in K,\|v\|_B\leq 1\big\}\text{ 
is bounded in $F$}\Bigr\}, \quad \text{ and } \\
C_b^{\{M\}}(U,F) &:= 
\Bigl\{f\in C^\infty(U,F):
\forall B\; \forall\text{ compact } K\subseteq U_B\;\exists \rho>0:\\
&\qquad\big\{\frac{f^{(k)}(x)(v_1,\dots,v_k)}{k!\,\rho^k\, M_k}:k\in \mathbb N,x\in K,\|v_i\|_B\leq 1\big\}\text{ 
is bounded in $F$}\Bigr\}
\\&
= \Bigl\{f\in C^\infty(U,F):
\forall B\; \forall\text{ compact } K\subseteq U_B\;\exists \rho>0:\\
&\qquad\big\{\frac{d^k_vf(x)}{k!\,\rho^k\, M_k}:k\in \mathbb N,x\in K,\|v\|_B\leq 1\big\}\text{ 
is bounded in $F$}\Bigr\}. 
\end{align*}
Here $B$ runs through all closed absolutely convex bounded subsets in $E$, $E_B$ is the vector space 
generated by $B$ with the Minkowski functional $\|v\|_B= \inf\{\lambda\ge 0: v\in \lambda B\}$ as complete norm, and 
$U_B = U \cap E_B$.
For Banach spaces $E$ and $F$ obviously
\[
C_b^{[M]}(U,F)=C^{[M]}(U,F).
\]
Now we define the spaces of main interest in this paper: 
\begin{equation*}\boxed{
\;C^{[M]}(U,F) := 
\Bigl\{f\in C^\infty(U,F):
\forall \ell \in F^*\; \forall B: \ell \o f \o i_B \in C^{[M]}(U_B,\mathbb{R})\Bigr\},
}\end{equation*}
where $B$ again runs through all closed absolutely convex bounded subsets in $E$, the mapping 
$i_B : E_B \to E$ denotes the inclusion
of $E_B$ in $E$, and $U_B = i_B^{-1}(U) = U \cap E_B$.
It will follow from Lemmas \ref{nmb:4.4} and \ref{nmb:4.5} that for Banach spaces 
$E$ and $F$ this definition coincides with the one given earlier in Section \ref{nmb:4.1}.

We equip $C^{[M]}(U,F)$ with the
initial locally convex structure induced by all linear mappings 
\begin{align*}
C^{[M]}(U,F) &\East{C^{[M]}(i_B,\ell)}{} C^{[M]}(U_B, \mathbb R), \quad f \mapsto \ell\o f\o i_B.
\end{align*}
Then $C^{[M]}(U,F)$ is a
convenient vector space as $c^\infty$-closed subspace in the product
$\prod_{\ell,B} C^{[M]}(U_B, \mathbb R)$,
since smoothness can be tested by composing with the inclusions $E_B\to E$ and with the
$\ell\in F^*$, see \cite[2.14.4 and 1.8]{KM97}.
This shows at the same time, that
\begin{align*}
C^{[M]}(U,F) &= 
\Bigl\{f\in F^U:
\forall \ell \in F^*\; \forall B: \ell \o f \o i_B \in C^{[M]}(U_B,\mathbb{R})\Bigr\}.
\end{align*}

\begin{lemma}[$C^{(M)}=C_b^{(M)}$]\label{nmb:4.4} 
Let $E,F$ be convenient vector spaces, and let $U \subseteq E$ be $c^\infty$-open.
Then a mapping $f : U \to F$ is $C^{(M)}$ (i.e., is in $C^{(M)}(U,F)$) if and
only if $f$ is $C_b^{(M)}$.
\end{lemma}

\begin{demo}{Proof}
Let $f : U \to F$ be $C^\infty$. We have the following equivalences, where $B$ runs through all closed absolutely convex bounded subsets in $E$:
\begin{align*}
f &\in C^{(M)}(U,F) \\
&\Longleftrightarrow 
\forall \ell \in F^* ~\forall B ~\forall K \subseteq U_B \text{ compact } ~\forall \rho>0:\\ 
&\hspace{1cm} \Big\{\frac{(\ell \o f)^{(k)}(x)(v_1,\ldots,v_k)}{\rho^k\, k!\, M_k} : x \in K, k \in \mathbb{N}, \|v_i\|_B\le 1\Big\} \text{ is bounded in } \mathbb{R}\\
&\Longleftrightarrow 
\forall B ~\forall K \subseteq U_B \text{ compact } ~\forall \rho>0 ~\forall \ell \in F^*:\\ 
&\hspace{1cm} \ell\Big(\Big\{\frac{f^{(k)}(x)(v_1,\ldots,v_k)}{\rho^k\, k!\, M_k} : x \in K, k \in \mathbb{N}, \|v_i\|_B\le 1\Big\}\Big) \text{ is bounded in } \mathbb{R}\\
&\Longleftrightarrow 
\forall B ~\forall K \subseteq U_B \text{ compact } ~\forall \rho>0: \\ 
&\hspace{1cm} \Big\{\frac{f^{(k)}(x)(v_1,\ldots,v_k)}{\rho^k\, k!\, M_k} : x \in K, k \in \mathbb{N}, \|v_i\|_B\le 1\Big\} \text{ is bounded in } F\\
&\Longleftrightarrow f \in C_b^{(M)}(U,F) \qed 
\end{align*}
\end{demo}

In the Roumieu case $C^{\{M\}}$ the corresponding equality holds only under additional assumptions:

\begin{lemma}[$C^{\{M\}}=C_b^{\{M\}}$]\label{nmb:4.5} 
Let $E,F$ be convenient vector spaces, and let $U \subseteq E$ be $c^\infty$-open.
Assume that there exists a Baire vector
space topology on the dual $F^*$ for which the point evaluations $\on{ev}_x$
are continuous for all $x\in F$.
Then a mapping $f : U \to F$ is $C^{\{M\}}$ if and
only if $f$ is $C_b^{\{M\}}$.
\end{lemma}

\begin{demo}{Proof}
$(\Rightarrow)$ 
Let $B$ be a closed absolutely convex bounded subset of $E$.
Let $K$ be compact in $U_B$.
We consider the sets
\[
A_{\rho,C} :=\Bigl\{\ell\in F^*: \frac{|(\ell\o f)^{(k)}(x)(v_1,\ldots,v_k)|}{\rho^{k}\,
k!\, M_{k}}\le C\text{ for all } x\in K, k\in \mathbb{N}, \|v_i\|_B\le 1\Bigr\}
\]
which are closed subsets in $F^*$ for the given Baire topology. We have
$\bigcup_{\rho,C}A_{\rho,C}= F^*$. By the Baire property there exist $\rho$
and $C$ such that the interior $\on{int}(A_{\rho,C})$ of $A_{\rho,C}$ is non-empty. If
$\ell_0\in \on{int}(A_{\rho,C})$, then for each $\ell\in F^*$ there is a $\delta>0$ such that 
$\delta\ell\in \on{int}(A_{\rho,C})-\ell_0$, and, hence, for all $k \in \mathbb{N}$, $x\in K$, and $\|v_i\|_B \le 1$, we have
\begin{align*}
|(\ell\o f)^{(k)}(x)(v_1,\ldots)| &\le \tfrac1\delta \Big(|((\delta\,\ell+\ell_0)\o f)^{(k)}(x)(v_1,\dots)| +
|(\ell_0\o f)^{(k)}(x)(v_1,\dots)|\Big)\\ &\le \tfrac{2C}{\delta}\,\rho^{k}\,k!\,
M_{k}.
\end{align*}
So the set 
\[
\Bigl\{\frac{f^{(k)}(x)(v_1,\ldots,v_k)}{\rho^{k}\, k!\, M_{k}}: x\in K, k\in\mathbb{N}, \|v_i\|_B\le 1\Bigr\}
\]
is weakly bounded in $F$ and hence bounded. Since $B$ and $K$ were arbitrary, we obtain $f \in C_b^{\{M\}}(U,F)$.

$(\Leftarrow)$ is obvious.
\qed\end{demo}

The following example shows that the additional assumption in Lemma \ref{nmb:4.5}
cannot be dropped.

\begin{example}\label{nmb:4.6}
By \cite[Thm.\ 1]{Thilliez08}, for each weakly log-convex sequence $M=(M_k)$ there exists 
$f \in C^{\{M\}}(\mathbb{R},\mathbb{R})$ such that $|f^{(k)}(0)| \ge k! \, M_k$ for all $k \in \mathbb{N}$. 
\emph{Then $g:\mathbb{R}^2 \to \mathbb{R}$ given by 
$g(s,t) = f(st)$ is $C^{\{M\}}$, 
whereas there is no reasonable
topology on $C^{\{M\}}(\mathbb{R},\mathbb{R})$ such that the associated mapping 
$g^\vee:\mathbb{R}\to C^{\{M\}}(\mathbb{R},\mathbb{R})$ is $C_b^{\{M\}}$.}
For a topology on $C^{\{M\}}(\mathbb{R},\mathbb{R})$ to be
reasonable we require only that all evaluations 
$\on{ev}_t: C^{\{M\}}(\mathbb{R},\mathbb{R}) \to \mathbb{R}$ are
bounded linear functionals.
\end{example}

\begin{demo}{Proof} 
The mapping $g$ is obviously $C^{\{M\}}$. If $g^\vee$ were $C_b^{\{M\}}$, for $s=0$ there existed $\rho$ such 
that 
\[
\Big\{\frac{(g^\vee)^{(k)}(0)}{k!\,\rho^k\,M_k}:k\in\mathbb N\Big\}
\]
was bounded in $C^{\{M\}}(\mathbb{R},\mathbb{R})$. We apply the bounded linear functional $\on{ev}_{t}$ for 
$t=2\rho$ and then get 
\[
\frac{|(g^\vee)^{(k)}(0)(2\rho)|}{k!\,\rho^k\,M_k}=\frac{(2\rho)^k |f^{(k)}(0)|}{k!\,\rho^k\,M_k} \ge 2^k,
\]
a contradiction.
\qed\end{demo}

This example shows that for $C_b^{\{M\}}$ one cannot expect cartesian closedness. 
Using cartesian closedness, i.e., Theorem \ref{nmb:5.2}, and Lemma \ref{nmb:5.1} this also shows (for 
$F=C^{\{M\}}(\mathbb{R},\mathbb{R})$ and $U=\mathbb{R}=E$) that
\[
C^{\{M\}}(U,F)\supsetneq\bigcap_{B,V} C_b^{\{M\}}(U_B,F_V),
\]
where $F_V$ is the completion 
of $F/p_V^{-1}(0)$ with respect to the seminorm $p_V$ induced by the absolutely convex closed 0-neighborhood 
$V$.

If we compose $g^\vee$ with the restriction mapping 
$(\on{incl}_\mathbb N)^*:C^{\{M\}}(\mathbb R,\mathbb R)\to\mathbb R^\mathbb N:=\prod_{t\in\mathbb N}\mathbb R$, then we get a $C^{\{M\}}$-curve, 
since the continuous linear functionals on $\mathbb R^\mathbb N$ are linear combinations of coordinate 
projections $\on{ev}_t$ with $t\in\mathbb N$. However, this curve cannot be $C_b^{\{M\}}$ as the argument above 
for $t>\rho$ shows.

\medskip

In the following Lemmas \ref{nmb:4.7} and \ref{nmb:4.8} we 
find projective descriptions for $C^{(M)}(U,F)$ and $C^{\{M\}}(U,F)$, if $E$, $F$ are Banach spaces, and $U \subseteq E$ is open.
This is of vital importance for the development of the convenient setting of $C^{\{M\}}(U,F)$. 
The spaces $C^{(M)}(U,F)$, however, already are projective by definition, and thus Lemma \ref{nmb:4.7} just 
gives a further projective description; see also Theorem \ref{nmb:8.6}. 
We include and use Lemma \ref{nmb:4.7} in order to treat the Beurling and Roumieu case in a uniform and efficient 
way.   

\begin{lemma} \label{nmb:4.7} 
Let $E,F$ be Banach spaces, $U \subseteq E$ open, and $f : U \to F$ a $C^\infty$-mapping.
The following are equivalent:
\begin{enumerate}
\item[(1)] $f$ is $C^{(M)}=C^{(M)}_b$.
\item[(2)] For each sequence $(r_k)$ with $r_k\, \rho^k \to 0$ for some $\rho>0$ and each compact $K \subseteq U$, 
the set 
\begin{equation*}
\Big\{\frac{f^{(k)}(a)(v_1,\ldots,v_k)}{k!M_k}\,r_k: a\in K, k\in \mathbb{N},\|v_i\| \le 1\Big\}
\end{equation*}
is bounded in $F$.
\item[(3)] For each sequence $(r_k)$ satisfying $r_k>0$, 
$r_kr_\ell\geq r_{k+\ell}$, 
and $r_k\, \rho^k \to 0$ for some $\rho>0$, each compact $K \subseteq U$, and each 
$\delta>0$, the set 
\begin{equation*}
\Big\{\frac{f^{(k)}(a)(v_1,\ldots,v_k)}{k!M_k}\,r_k\, \delta^k: a\in K, k\in \mathbb{N},\|v_i\| \le 1\Big\}
\end{equation*}
is bounded in $F$. 
\end{enumerate}
\end{lemma}

\begin{demo}{Proof}
\thetag{1} $\Rightarrow$ \thetag{2}
For $(r_k)$ and $K$, and $\rho>0$ such that $r_k\, \rho^k \to 0$,
\[
\left\|\frac{f^{(k)}(a)}{k!\,M_k}r_k \right\|_{L^k(E,F)} = 
\left\|\frac{f^{(k)}(a)}{k!\,\rho^k\,M_k}\right\|_{L^k(E,F)} \cdot |r_k \rho^k| 
\]
is bounded uniformly in $k\in \mathbb N$ and $a\in K$ (by Lemma \ref{nmb:4.4}). 

\thetag{2} $\Rightarrow$ \thetag{3} Apply \thetag{2} to the sequence $(r_k \delta^k)$.

\thetag{3} $\Rightarrow$ \thetag{1} 
Let $a_k:=\sup_{a\in K}\|\frac{f^{(k)}(a)}{k!\,M_k}\|_{L^k(E,F)}$. 
By the following lemma, the $a_k$ are the coefficients of a power series with infinite radius of convergence.
Thus $a_k/\rho^k$ is bounded for every $\rho>0$. 
\qed\end{demo}

\begin{lemma*}
For a formal power series $\sum_{k \ge 0} a_k t^k$ with real coefficients the following are equivalent:
\begin{enumerate}
\item[\thetag{4}] The radius of convergence is infinite.
\item[\thetag{5}] For each sequence $(r_k)$ satisfying $r_k>0$, 
$r_kr_\ell\geq r_{k+\ell}$, 
and $r_k\,\rho^k \to 0$ for some $\rho>0$, and each 
$\delta>0$, the sequence $(a_k r_k \delta^k)$ is bounded.
\end{enumerate}
\end{lemma*}

\begin{demo}{Proof}
\thetag{4} $\Rightarrow$ \thetag{5} The series 
$\sum a_k r_k \delta^k = \sum (a_k (\frac{\delta}{\rho})^k) r_k \rho^k$ converges absolutely for each $\delta$. 
Hence $(a_k r_k \delta^k)$ is bounded.

\thetag{5} $\Rightarrow$ \thetag{4} Suppose that the radius of convergence $\rho$ is finite. 
So $\sum_k |a_k| n^k = \infty$ for $n > \rho$.
Set $r_k = 1/n^k$. Then, by \thetag{5}, 
\[
a_k n^k 2^k = a_k r_k n^{2k} 2^k = a_k r_k (2 n^2)^k < C,
\] 
for some $C>0$ and all $k$. Consequently, $\sum_k |a_k| n^k \le C \sum_k \frac{1}{2^{k}}$, a contradiction.
\qed\end{demo}

\begin{lemma} \label{nmb:4.8} 
Let $E,F$ be Banach spaces, $U \subseteq E$ open, and $f : U \to F$ a $C^\infty$-mapping.
The following are equivalent:
\begin{enumerate}
\item[(1)] $f$ is $C^{\{M\}}=C_b^{\{M\}}$.
\item[(2)] For each sequence $(r_k)$ with $r_k\,\rho^k\to
0$ for all $\rho>0$, and each compact $K\subseteq U$, 
the set 
\[
\Big\{\frac{f^{(k)}(a)(v_1,\dots,v_k)}{k!M_k}\,r_k: a\in K, k\in \mathbb{N},\|v_i\| \le 1\Big\}
\] 
is bounded in $F$.
\item[(3)] For each sequence $(r_k)$ satisfying $r_k>0$, 
$r_kr_\ell\geq r_{k+\ell}$, 
and $r_k\,\rho^k\to 0$ for all $\rho>0$, and each compact $K \subseteq U$, 
there exists 
$\delta>0$ such that 
\[
\Big\{\frac{f^{(k)}(a)(v_1,\ldots,v_k)}{k!M_k}\,r_k\,\delta^k: a\in K, k\in \mathbb{N},\|v_i\| \le 1\Big\}
\] 
is bounded in $F$. 
\end{enumerate}
\end{lemma}

\demo{Proof}
\thetag{1} $\Rightarrow$ \thetag{2}
For $K$, there exists $\rho>0$ such that
\[
\Big\|\frac{f^{(k)}(a)}{k!\,M_k}r_k \Big\|_{L^k(E,F)} = 
\Big\|\frac{f^{(k)}(a)}{k!\,\rho^k\,M_k}\Big\|_{L^k(E,F)} \cdot |r_k\rho^k| 
\]
is bounded uniformly in $k\in \mathbb N$ and $a\in K$ (by Lemma \ref{nmb:4.5}). 

\thetag{2} $\Rightarrow$ \thetag{3} Use $\delta=1$.

\thetag{3} $\Rightarrow$ \thetag{1} 
Let $a_k:=\sup_{a\in K}\|\frac{f^{(k)}(a)}{k!\,M_k}\|_{L^k(E,F)}$. Using 
\cite[9.2(4$\Rightarrow$1)]{KM97} 
these are the coefficients of a power series with positive 
radius of convergence.
Thus $a_k/\rho^k$ is bounded for some $\rho>0$. 
\qed\enddemo

\begin{definition}[Banach plots]\label{nmb:4.9}
Let $E$ be a convenient vector space. 
A \emph{$C^{[M]}$ (Banach) plot in $E$} is a mapping $c : D \to E$ of class $C^{[M]}$, 
where $D$ is an open set in some Banach space $F$. It suffices to only consider the open unit ball 
$D=oF$.
\end{definition}

\begin{theorem} \label{nmb:4.10}
Let $M=(M_k)$ be log-convex.
Let $U\subseteq E$ be $c^\infty$-open in a convenient vector space $E$, let $F$
be a Banach space, and let $f:U\to F$ be a mapping. Then:
\begin{align*}
f \in C^{[M]}(U,F) &\implies f \o c \in C^{[M]}, \text{ for all $C^{[M]}$-plots } c.
\end{align*} 
\end{theorem}
Note that the converse ($\Leftarrow$) holds by Section \ref{nmb:4.3}.

\begin{demo}{Proof}
We treat first the Beurling case $C^{(M)}$: 
We have to show that $f\o c$ is $C^{(M)}$ for each $C^{(M)}$-plot $c: G \supseteq D\to E$. 
By Lemma \ref{nmb:4.7}(3), it suffices to show that, for each 
sequence $(r_k)$ satisfying $r_k>0$, 
$r_kr_\ell\geq r_{k+\ell}$, 
and $r_k\, t^k\to 0$ for some $t>0$, each compact $K \subseteq D$, 
and each
$\delta>0$, the set
\begin{equation} \label{eq:13}
\Big\{\frac{(f\o c)^{(k)}(a)(v_1,\ldots,v_k)}{k!M_k}\,r_k\,\delta^k: a\in K, k\in \mathbb{N}, \|v_i\|_G \le 1\Big\} 
\end{equation}
is bounded in $F$.

So let $\delta$, the sequence $(r_k)$, and a compact (and without loss of generality convex) subset $K\subseteq D$ be fixed.
For each $\ell\in E^*$ the set 
\begin{equation} \label{eq:14}
\Big\{\frac{(\ell\o c)^{(k)}(a)(v_1,\ldots,v_k)}{k!M_k}\,r_k\,(2\delta)^k: a\in K, k\in \mathbb{N},\|v_i\|_G \le 1\Big\}
\end{equation}
is bounded in $\mathbb{R}$, by Lemma \ref{nmb:4.7}(2) applied to the sequence $(r_k(2\delta)^k)$.
Thus, the set
\begin{equation} \label{eq:15}
\Big\{\frac{c^{(k)}(a)(v_1,\ldots,v_k)}{k!M_k}\,r_k\,(2\delta)^k: a\in K, k\in \mathbb{N},\|v_i\|_G \le 1\Big\}
\end{equation}
is contained in some closed absolutely convex bounded subset $B$ of $E$ and hence
\[
\frac{\|c^{(k)}(a)\|_{L^k(G,E_B)} \, r_k\,\delta^k}{k!M_k}
\leq \frac1{2^k}. 
\]
Furthermore $c:K\to E_B$ is Lipschitzian, since
\[
c(x)-c(y)=\int_0^1 c'(y+t(x-y))\,(x-y)\,dt\in \frac{M_1\,\|x-y\|_G}{2\, r_1\,\delta}\,B, 
\]
and hence $c(K)$ is compact in $E_B$.
By Fa\`a di Bruno's formula for Banach spaces (see \cite{FaadiBruno1855} for the 1-dimensional version), for $k\ge 1$,
\begin{align*}
\frac{(f\o c)^{(k)}(a)}{k!} = \on{sym}\Big(
\sum_{j\ge 1} \sum_{\substack{\alpha\in \mathbb{N}_{>0}^j\\ \alpha_1+\dots+\alpha_j =k}}
\frac{1}{j!}f^{(j)}(c(a))\Big( 
\frac{c^{(\alpha_1)}(a)}{\alpha_1!},\dots,
\frac{c^{(\alpha_j)}(a)}{\alpha_j!}\Big)\Big),
\end{align*}
where $\on{sym}$ denotes symmetrization.
Using \eqref{eq:4} and $f \in C^{(M)}(U,F)$, we find that for each $\rho>0$ there is $C>0$ so that,
for all $a\in K$ and $k\in \mathbb{N}_{>0}$, 
\begin{align} 
\Big\|&\frac{(f\o c)^{(k)}(a)}{k!M_k}\, r_k\,\delta^k \Big\|_{L^k(G,F)} \notag 
\\
&\le \sum_{j\ge 1} M_1^j \!\!\!\!\sum_{\substack{\alpha\in \mathbb{N}_{>0}^j\\ \alpha_1+\dots+\alpha_j =k}}
\underbrace{\frac{\|f^{(j)}(c(a))\|_{L^j(E_B,F)}}{j!M_j}}_{\leq C\,\rho^j}\;\prod_{i=1}^j\;
\underbrace{\frac{\|c^{(\alpha_i)}(a)\|_{L^{\alpha_i}(G,E_B)} \, 
r_{\alpha_i}\,\delta^{\alpha_i}}{\alpha_i!M_{\alpha_i}}}_{\leq 1/2^{\alpha_i}} \notag
\\ 
&\le \sum_{j\ge 1} M_1^j \binom{k-1}{j-1} C\,\rho^j\, \frac1{2^k}
= M_1 \rho(1+M_1\,\rho)^{k-1} C\, \frac1{2^k}\leq \frac{C}2, \label{eq:16}
\end{align}
as required, where in the last inequality we set 
$\rho:=1/{M_1}$.

Let us now consider the Roumieu case $C^{\{M\}}$:
Let now $c: G \supseteq D\to E$ be a $C^{\{M\}}$-plot.
We have to show that $f\o c$ is $C^{\{M\}}$. 
By Lemma \ref{nmb:4.8}(3), it suffices to show that for each 
sequence $(r_k)$ satisfying $r_k>0$, 
$r_kr_\ell\geq r_{k+\ell}$, 
and $r_k\,t^k\to 0$ for all $t>0$, and each compact $K \subseteq D$, 
there exists $\delta>0$ such that the set \eqref{eq:13} is bounded in $F$.

By Lemma \ref{nmb:4.8}(2) (applied to $(r_k2^k)$ instead of $(r_k)$), for each $\ell\in E^*$, 
each sequence $(r_k)$ with $r_k\,t^k\to 0$ for all $t>0$, and each compact $K \subseteq D$, the set 
\eqref{eq:14} with $\delta=1$
is bounded in $\mathbb R$, and,
thus, the set \eqref{eq:15} with $\delta=1$ 
is contained in some closed absolutely convex bounded subset $B$ of $E$.
Using that $f \in C^{\{M\}}(U,F)$ and computing as above we find that, for some $\rho>0$ and $C>0$ and $\delta:=\frac{2}{1+M_1\,\rho}$, 
the left-hand side of 
\eqref{eq:16}
is bounded by $\frac{C\, M_1\, \rho}{1+M_1\,\rho}$.
\qed\end{demo}

\begin{theorem}[$C^{[M]}$ is a category]\label{nmb:4.11}
Let $M=(M_k)$ be log-convex.
Let $E,F,G$ be convenient vector spaces, $U \subseteq E$, $V \subseteq F$ be $c^\infty$-open,
and $f : U \to F$, $g : V \to G$, and $f(U) \subseteq V$. Then:
\begin{align*}
f,g \in C^{[M]} &\implies g\o f \in C^{[M]}.
\end{align*}
\end{theorem}

\begin{demo}{Proof}
By Section \ref{nmb:4.3}, we must show that for all closed absolutely convex bounded $B \subseteq E$ 
and for all $\ell \in G^*$ the composite
$\ell \o g \o f \o i_B : U_B \to \mathbb{R}$ belongs to $C^{[M]}$.
\[
\xymatrix{
U \ar[rr]^{f} && V \ar[rr]^g \ar[drr]_{\ell \o g} && G \ar[d]^{\ell} \\
U_B \ar[u]^{i_B} \ar[urr]_{f \o i_B} &&&& \mathbb{R}
}
\]
By assumption, $f \o i_B$ and $\ell \o g$ are $C^{[M]}$.
So the assertion follows from Theorem \ref{nmb:4.10}. 
\qed\end{demo}

\section{\label{nmb:5}The exponential law}

\begin{lemma} \label{nmb:5.1}
Let $E$ be a Banach space, and $U \subseteq E$ be open.
Let $F$ be a convenient vector space, and let $\mathcal S$ be a
family of bounded linear functionals on $F$ which together detect bounded
sets (i.e., $B\subseteq F$ is bounded if and only if $\ell(B)$ is bounded for all
$\ell\in\mathcal S$). Then:
\begin{align*}
f \in C^{[M]}(U,F) &\Longleftrightarrow \ell\o f \in C^{[M]}(U,\mathbb{R}), \text{ for all } 
\ell \in \mathcal{S}.
\end{align*}
\end{lemma}

\begin{demo}{Proof}
For $C^\infty$-curves this follows from \cite[2.1 and 2.11]{KM97}, and, by composing with such, 
it follows for $C^\infty$-mappings $f: U \to F$.

In the Beurling case $C^{(M)}$:
By Lemma \ref{nmb:4.7}, for $\ell\in F^*$, 
the function $\ell\o f$ is $C^{(M)}$ if and only if 
the set 
\begin{equation} \label{eq:17}
\Big\{\frac{(\ell\o f)^{(k)}(a)(v_1,\dots,v_k)}{k!M_k}\,r_k : a\in K, k\in \mathbb{N},\|v_i\|\le 1\Big\} 
\end{equation}
is bounded, for each sequence $(r_k)$ with $r_k\, \rho^k \to 0$ for some $\rho >0$ and 
each compact $K \subseteq U$.
So the smooth mapping $f: U \to F$ is $C^{(M)}$ if and only if 
the set 
\begin{equation} \label{eq:18}
\Big\{\frac{f^{(k)}(a)(v_1,\dots,v_k)}{k!M_k} \,r_k : a\in K, k\in \mathbb{N},\|v_i\|\le 1\Big\} 
\end{equation}
is bounded in $F$, for each such $(r_k)$ and $K$.
This is in turn equivalent to $\ell\o f\in C^{(M)}$ for all $\ell\in\mathcal S$, 
since $\mathcal S$ detects bounded sets.

The same proof works in the Roumieu case $C^{\{M\}}$ if we use Lemma \ref{nmb:4.8} and 
demand that $r_k\,\rho^k\to 0$ for all $\rho>0$.
\qed\end{demo}

\begin{theorem}[Cartesian closedness]\label{nmb:5.2} 
We have: 
\begin{enumerate}
\item[\thetag{1}] Let $M=(M_k)$ be weakly log-convex and have moderate growth.
Then, for convenient vector spaces $E_1$, $E_2$, and $F$
and $c^\infty$-open sets $U_1\subseteq E_1$ and $U_2\subseteq E_2$, we have the exponential law:
\begin{align*}
f \in C^{[M]}(U_1 \times U_2,F) &\Longleftrightarrow f^\vee \in C^{[M]}(U_1,C^{[M]}(U_2,F)).
\end{align*}
The direction \thetag{$\Leftarrow$} holds 
without the assumption that $M=(M_k)$ has moderate growth. 
The direction \thetag{$\Rightarrow$} holds 
without the assumption that $M=(M_k)$ is weakly log-convex.
\item[\thetag{2}] Let $M=(M_k)$ be log-convex and have moderate growth.
Then the category of $C^{[M]}$-mappings between convenient real vector spaces 
is cartesian closed, i.e., satisfies the exponential law.
\end{enumerate} 
\end{theorem}

Note that $C^{[M]}$ is not necessarily a category if $M=(M_k)$ is just weakly log-convex. 

\begin{demo}{Proof}
\thetag{2} is a direct consequence of \thetag{1} and Theorem \ref{nmb:4.11}. Let us 
prove \thetag{1}. 
We have $C^{\infty}(U_1 \times U_2,F) \cong C^\infty(U_1,C^\infty(U_2,F))$, by \cite[3.12]{KM97}; 
thus, in the following all mappings are assumed to be smooth. 
We have the following equivalences, where $B \subseteq E_1 \times E_2$ 
and $B_i \subseteq E_i$ run through all closed absolutely convex bounded subsets, respectively: 
\begin{align*}
f &\in C^{[M]}(U_1 \times U_2,F) \\
&\Longleftrightarrow 
\forall \ell \in F^* ~\forall B: \ell \o f \o i_B \in C^{[M]}((U_1 \times U_2)_B,\mathbb{R})\\
&\Longleftrightarrow 
\forall \ell \in F^* ~\forall B_1,B_2: \ell \o f \o (i_{B_1} \times i_{B_2}) \in C^{[M]}((U_1)_{B_1} \times (U_2)_{B_2},\mathbb{R}) 
\end{align*}
For the second equivalence we use that every bounded $B \subseteq E_1 \times E_2$ is contained in $B_1 \times B_2$ for some bounded $B_i \subseteq E_i$, 
and, thus, the inclusion $(E_1 \times E_2)_B \to (E_1)_{B_1} \times (E_1)_{B_2}$ is bounded.

On the other hand, we have:
\begin{align*}
f^\vee &\in C^{[M]}(U_1,C^{[M]}(U_2,F)) \\
&\Longleftrightarrow 
\forall B_1: f^\vee \o i_{B_1} \in C^{[M]}((U_1)_{B_1},C^{[M]}(U_2,F))\\
&\Longleftrightarrow 
\forall \ell \in F^* ~\forall B_1,B_2: C^{[M]}(i_{B_2},\ell) \o f^{\vee} \o i_{B_1} \in C^{[M]}((U_1)_{B_1},C^{[M]}((U_2)_{B_2},\mathbb{R})) 
\end{align*} 
For the second equivalence we use Lemma \ref{nmb:5.1} and the fact that the linear mappings $C^{[M]}(i_{B_2},\ell)$ generate the bornology. %

These considerations imply that in order to prove cartesian closedness we may restrict to the case that 
$U_i \subseteq E_i$ are open in Banach spaces $E_i$ and $F = \mathbb{R}$.

\medskip

({\bf Direction} $\Rightarrow$) We assume that $M=(M_k)$ has moderate growth.
Let $f \in C^{[M]}(U_1 \times U_2,\mathbb{R})$. It is clear that $f^\vee$ takes values in $C^{[M]}(U_2,\mathbb{R})$.

\medskip

\noindent
{\bf Claim.} {\it $f^\vee:U_1\to C^{[M]}(U_2,\mathbb R)$ is $C^\infty$ with $d^jf^\vee=(\partial_1^jf)^\vee$.}
\newline
Since $C^{[M]}(U_2,\mathbb R)$ is a convenient vector space, by \cite[5.20]{KM97}
it is enough to show that the iterated unidirectional
derivatives $d^j_vf^\vee(x)$ exist, equal $\partial_1^jf(x,\quad)(v^j)$, 
and are separately bounded for $x$, resp.\  $v$, in compact 
subsets.
For $j=1$ and fixed $x$, $v$, and $y$ consider the smooth curve
$c:t\mapsto f(x+tv,y)$. By the fundamental theorem
\begin{align*}
\frac{f^\vee(x+tv)-f^\vee(x)}{t}(y)&-(\partial_1f)^\vee(x)(y)(v) 
= \frac{c(t)-c(0)}t-c'(0) \\
&= t\int_0^1 s \int_0^1 c''(tsr)\,dr\,ds \\
&= t\int_0^1 s\int_0^1 \partial_1^2f(x+tsrv,y)(v,v)\,dr \,ds.
\end{align*}
Since $(\partial_1^2f)^\vee(K_1)(oE_1^2)$ is obviously bounded in $C^{[M]}(U_2,\mathbb R)$ for each compact subset
$K_1\subseteq U_1$, this expression is Mackey convergent 
to 0 in $C^{[M]}(U_2,\mathbb R)$ as $t\to 0$.
Thus $d_vf^\vee(x)$ exists and equals $\partial_1f(x,\quad)(v)$.

Now we proceed by induction, applying the same arguments as before to
$(d^j_vf^\vee)^\wedge : (x,y)\mapsto \partial_1^jf(x,y)(v^j)$ instead of $f$.
Again
$(\partial_1^2(d^j_vf^\vee)^\wedge)^\vee(K_1)(oE_1^2) = (\partial_1^{j+2}f)^\vee(K_1)(oE_1,oE_1,v,\dots,v)$
is bounded, and also the separated boundedness of $d^j_vf^\vee(x)$ follows.
So the claim is proved.
\medskip

We have to show that $f^\vee : U_1 \to C^{[M]}(U_2,\mathbb{R})$ is $C^{[M]}$.

In the Beurling case $C^{(M)}$: 
\begin{equation} \label{eq:19}
\begin{split}
\xymatrix{
U_1\ar@{->}[r]^(0.3){f^\vee} &C^{(M)}(U_2,\mathbb{R} )\ar@{=}[r] &
\varprojlim _{K_2}\varprojlim _{\rho _2}C^{M}_{\rho _2}(E_2\supseteq K_2,\mathbb{R} )
\ar@{->}[r]^(0.8){\ell} \ar@{->}[d] &\mathbb{R} \\
& &\varprojlim _{\rho _2}C^{M}_{\rho _2}(E_2\supseteq K_2,\mathbb{R} )\ar@{->}[d] &
\\
K_1\ar@{^{ (}->}[uu] \ar@{.>}[r] &C^M_{K_2,\rho_2}(U_2,\mathbb R)
\ar@{->}[r] &C^{M}_{\rho _2}(E_2\supseteq K_2,\mathbb{R} )
\ar@/_2pc/@{.>}[uur] & \\
} 
\end{split} 
\end{equation}
By Lemma \ref{nmb:5.1}, it suffices to show that $f^\vee : U_1 \to C^{M}_{\rho_2}(E_2\supseteq K_2,\mathbb{R})$ is 
$C_b^{(M)}= C^{(M)}$ (see Lemma \ref{nmb:4.4}) for each compact $K_2 \subseteq U_2$ and each $\rho_2>0$, 
since every $\ell \in C^{(M)}(U_2,\mathbb{R})^*$ factors over some $C^{M}_{\rho_2}(E_2\supseteq K_2,\mathbb{R})$.
Thus it suffices to prove that, for all compact $K_1 \subseteq U_1$, $K_2 \subseteq U_2$ and all $\rho_1,\rho_2>0$, 
the set 
\begin{equation} \label{eq:20}
\Big\{\frac{d^{k_1} f^\vee(x_1)(v_1^1,\dots,v_{k_1}^1)}{k_1!\, \rho_1^{k_1}\, M_{k_1}} : x_1 \in K_1, k_1 \in \mathbb{N}, \|v_j^1\|_{E_1} \le 1\Big\}
\end{equation}
is bounded in $C^{M}_{K_2,\rho_2}(U_2,\mathbb{R})$, or, equivalently, for all $K_1,K_2,\rho_1,\rho_2$
the set 
\begin{equation} \label{eq:21}
\Big\{\frac{\partial_2^{k_2}\partial_1^{k_1} f(x_1,x_2)(v_1^1,\dots,v_{k_1}^1;v_1^2,\dots,v_{k_2}^2)}{k_2!\,k_1!\, \rho_2^{k_2}\,\rho_1^{k_1}\, M_{k_2}\, M_{k_1}} : x_i \in K_i, k_i \in \mathbb{N}, \|v_j^i\|_{E_i} \le 1\Big\}
\end{equation}
is bounded in $\mathbb{R}$.

Since $M=(M_k)$ has moderate growth, i.e., $M_{k_1+k_2} \le \sigma^{k_1+k_2} M_{k_1} M_{k_2}$ for some $\sigma>0$, we 
obtain, for $x_1 \in K_1$, $k_1 \in \mathbb{N}$, and $\|v_j^1\|_{E_1}\le 1$, 
\begin{align} \label{eq:22}
\begin{split}
	&\Big\|\frac{d^{k_1} f^\vee(x_1)(v_1^1,\dots,v_{k_1}^1)}{k_1!\, \rho_1^{k_1}\, M_{k_1}}\Big\|_{K_2,\rho_2} 
\\
&= \sup \Big\{\frac{|\partial_2^{k_2}\partial_1^{k_1} f(x_1,x_2)(v_1^1,\dots,v_{k_1}^1;v_1^2,\dots,v_{k_2}^2)|}
{k_2!\,k_1!\, \rho_2^{k_2}\,\rho_1^{k_1}\, M_{k_2}\, M_{k_1}} : 
x_2 \in K_2, k_2 \in \mathbb{N}, \|v_j^2\|_{E_2} \le 1\Big\}  %\notag
\\
&\le \sup \Big\{(2\sigma)^{k_1+k_2}\frac{|\partial_2^{k_2}\partial_1^{k_1} f(x_1,x_2)(v_1^1,\dots;v_1^2,\dots)|}
{(k_1+k_2)!\, \rho_1^{k_1}\,\rho_2^{k_2}\, M_{k_1+k_2}} : 
x_2 \in K_2, k_2 \in \mathbb{N}, \|v_j^2\|_{E_2} \le 1\Big\}.  %\notag
\end{split}
\end{align}
If for given $\rho_1,\rho_2>0$ we set $\rho := \tfrac{1}{2 \sigma} \min \{\rho_1,\rho_2\}$, then 
\eqref{eq:22} is bounded by
\begin{equation} \label{eq:23}
\sup \Big\{\frac{|\partial_2^{k_2}\partial_1^{k_1} f(x_1,x_2)(v_1^1,\dots;v_1^2,\dots)|}{(k_1+k_2)!\, \rho^{k_1+k_2}\, M_{k_1+k_2}} : 
x_i \in K_i, k_i \in \mathbb{N}, \|v_j^i\|_{E_i} \le 1\Big\}
\end{equation}
which is finite, since $f$ is $C^{(M)}$. Thus, $f^\vee$ is $C^{(M)}$. 

In the Roumieu case $C^{\{M\}}$: 
\begin{equation} \label{eq:24}
\begin{split}
\xymatrix{
U_1\ar@{->}[r]^(0.3){f^\vee} &C^{\{M \} }(U_2,\mathbb{R} )\ar@{=}[r] &
\varprojlim _{K_2}\varinjlim _{\rho _2}C^{M}_{\rho _2}(E_2\supseteq K_2,\mathbb{R} )
\ar@{->}[r]^(0.7){\ell} \ar@{->}[d] &\mathbb{R} \\
& &\varinjlim _{\rho _2}C^{M}_{\rho _2}(E_2\supseteq K_2,\mathbb{R} )\ar@{.>}[ur] &
\\
K_1\ar@{^{ (}->}[uu] \ar@{.>}[r] &C^M_{K_2,\rho_2}(U_2,\mathbb R)\ar@{->}[r] &
C^{M}_{\rho _2}(E_2\supseteq K_2,\mathbb{R} )\ar@{->}[u] & \\
}
\end{split} 
\end{equation}
By Lemma \ref{nmb:5.1}, it suffices to show that $f^\vee : U_1 \to \varinjlim_{\rho_2} 
C^{M}_{\rho_2}(E_2\supseteq K_2,\mathbb{R})$ is $C_b^{\{M\}}\subseteq C^{\{M\}}$ for each compact $K_2 \subseteq U_2$, 
since every $\ell \in C^{\{M\}}(U_2,\mathbb{R})^*$ factors over some $\varinjlim_{\rho_2} C^{M}_{\rho_2}(E_2\supseteq K_2,\mathbb{R})$.
Thus it suffices to prove that, for all compact $K_1 \subseteq U_1$ and $K_2 \subseteq U_2$ there exists $\rho_1>0$, 
such that the set \eqref{eq:20} 
is bounded in $\varinjlim_{\rho_2} C^{M}_{\rho_2}(E_2\supseteq K_2,\mathbb{R})$. For that it suffices to show that 
for all $K_1,K_2$ there are $\rho_1,\rho_2$ so that  
the set \eqref{eq:21}
is bounded in $\mathbb{R}$.

Since $f$ is $C^{\{M\}}$, there exists $\rho>0$ so that \eqref{eq:23} is finite, by  Proposition  
\ref{nmb:4.2}(3). If we set $\rho_i := 2 \sigma \rho$, then 
\eqref{eq:22} is bounded by \eqref{eq:23}. It follows that $f^\vee$ is $C^{\{M\}}$.

\medskip

({\bf Direction} $\Leftarrow$) Let $f^\vee : U_1 \to C^{[M]}(U_2,\mathbb{R})$ be $C^{[M]}$. 
Clearly, $f^\vee:U_1\to C^{[M]}(U_2,\mathbb{R})\to C^\infty(U_2,\mathbb{R})$
is $C^\infty$, see Proposition \ref{nmb:8.1}, and so it remains to show that $f \in C^{[M]}(U_1 \times U_2,\mathbb{R})$.

In the Beurling case $C^{(M)}$: Consider diagram \eqref{eq:19}.
For each compact $K_2 \subseteq U_2$ and each $\rho_2>0$, the mapping $f^\vee : U_1 
\to C^{M}_{\rho_2}(E_2\supseteq K_2,\mathbb{R})$ is $C^{(M)}=C_b^{(M)}$.
That means that, for all compact $K_1 \subseteq U_1$, $K_2 \subseteq U_2$ and all $\rho_1,\rho_2>0$,
the set \eqref{eq:20} is bounded in $C^{M}_{\rho_2}(E_2\supseteq K_2,\mathbb{R})$.
Since it is contained
in $C^{M}_{K_2,\rho_2}(U_2,\mathbb R)=\{g\in C^\infty(U_2,\mathbb R):j^\infty g|_{K_2}\in C^{M}_{\rho_2}(E_2\supseteq K_2,\mathbb{R})\}$
and $\|g\|_{K_2,\rho_2}=\|j^\infty g|_{K_2}\|_{\rho_2}$, it is also bounded in this space, and hence
the set \eqref{eq:21} is bounded.

Since $M=(M_k)$ is weakly log-convex, thus, $k_1!\,k_2!\, M_{k_1} M_{k_2} \le (k_1+k_2)!\, M_{k_1+k_2}$, we 
have, for $x_1 \in K_1$, $k_1 \in \mathbb{N}$, and $\|v_j^1\|_{E_1}\le 1$, 
\begin{align} \label{eq:25}
\begin{split}
	&\Big\|\frac{d^{k_1} f^\vee(x_1)(v_1^1,\dots,v_{k_1}^1)}{k_1!\, \rho_1^{k_1}\, M_{k_1}}\Big\|_{K_2,\rho_2} 
\\
&= \sup \Big\{\frac{|\partial_2^{k_2}\partial_1^{k_1} f(x_1,x_2)(v_1^1,\dots,v_{k_1}^1;v_1^2,\dots,v_{k_2}^2)|}
{k_2!\,k_1!\, \rho_2^{k_2}\,\rho_1^{k_1}\, M_{k_2}\, M_{k_1}} : 
x_2 \in K_2, k_2 \in \mathbb{N}, \|v_j^2\|_{E_2} \le 1\Big\} %\notag
\\
&\ge \sup \Big\{\frac{|\partial_2^{k_2}\partial_1^{k_1} f(x_1,x_2)(v_1^1,\dots;v_1^2,\dots)|}
{(k_1+k_2)!\, \rho_1^{k_1}\,\rho_2^{k_2}\, M_{k_1+k_2}} : 
x_2 \in K_2, k_2 \in \mathbb{N}, \|v_j^2\|_{E_2} \le 1\Big\}. %\notag
\end{split}
\end{align}
This implies that $f$ is $C^{(M)}$.

In the Roumieu case $C^{\{M\}}$: Consider diagram \eqref{eq:24}.
For each compact $K_2 \subseteq U_2$, the mapping $f^\vee : U_1 \to \varinjlim_{\rho_2} C^{M}_{\rho_2}(E_2\supseteq K_2,\mathbb{R})$ is $C^{\{M\}}$.
The inductive limit is regular, by Proposition \ref{nmb:4.2}(3).
So the dual space $(\varinjlim_{\rho_2} C^{M}_{\rho_2}(E_2\supseteq K_2,\mathbb{R}))^*$ 
can be equipped with the Baire topology of the countable limit 
$\varprojlim_{\rho_2} C^{M}_{\rho_2}(E_2\supseteq K_2,\mathbb{R})^*$ of Banach spaces. 
Thus, the mapping $f^\vee : U_1 \to \varinjlim_{\rho_2} C^{M}_{\rho_2}(E_2\supseteq K_2,\mathbb{R})$ is $C_b^{\{M\}}$, by Lemma \ref{nmb:4.5}.
By regularity, for each compact $K_1 \subseteq U_1$ there exists $\rho_1 >0$ so that the set \eqref{eq:20} 
is contained and bounded in
$C^{M}_{\rho_2}(E_2\supseteq K_2,\mathbb{R})$ for some $\rho_2>0$.
Since this set is contained in 
$C^{M}_{K_2,\rho_2}(U_2,\mathbb R)=\{g\in C^\infty(U_2,\mathbb R):j^\infty g|_{K_2}\in C^{M}_{\rho_2}(E_2\supseteq K_2,\mathbb{R})\}$
and $\|g\|_{K_2,\rho_2}=\|j^\infty g|_{K_2}\|_{\rho_2}$, it is also bounded in this space, and hence
the set \eqref{eq:21} is bounded.
Then \eqref{eq:25} implies that $f$ is $C^{\{M\}}$. The proof is complete.
\qed\end{demo}

\begin{remarks} \label{nmb:5.3}
Theorem \ref{nmb:8.2} below states that, if $M=(M_k)$ is (weakly) log-convex, $E, F$ are convenient vector spaces, 
and $U \subseteq E$ is $c^\infty$-open, then 
\begin{equation} \label{eq:26}
C^{\{M\}}(U,F) = \varprojlim_{L} C^{(L)}(U,F)
\end{equation}
as vector spaces with bornology,
where the projective limits are taken over all (weakly) log-convex $L=(L_k)$ with $M \lhd L$.
Using this equality we can give an alternative proof of the direction
\[
f^\vee \in C^{\{M\}}(U_1, C^{\{M\}}(U_2,F)) \Rightarrow f \in C^{\{M\}}(U_1 \times U_2, F)
\]
in Theorem \ref{nmb:5.2} as follows: 
If $f^\vee \in C^{\{M\}}(U_1, C^{\{M\}}(U_2,F))$ then $f^\vee \in C^{(L)}(U_1, C^{(L)}(U_2,F))$ for all $L=(L_k)$ with 
$M \lhd L$, by \eqref{eq:26}. By cartesian closedness of $C^{(L)}$ 
(i.e., Theorem \ref{nmb:5.2}, the implication which holds without moderate growth), we 
have $f \in C^{(L)}(U_1 \times U_2, F)$ for all $L$, and, by \eqref{eq:26} again, 
$f \in C^{\{M\}}(U_1 \times U_2, F)$.

The proof of \eqref{eq:26} in Theorem \ref{nmb:8.2} uses the 
$C^{\{M\}}$ uniform boundedness principle, i.e., Theorem \ref{nmb:6.1}, and the proof of the latter 
uses completeness of the inductive limit $\varinjlim_{\rho} C^M_{\rho}(E \supseteq K,F)$, where $E, F$ are Banach spaces and 
$K \subseteq E$ is compact, see Proposition \ref{nmb:4.2}.
Here is a direct proof of \eqref{eq:26}, where we only assume that $M=(M_k)$ is positive: 

The spaces coincide as vector spaces by Sections \ref{nmb:4.1}, \ref{nmb:4.3}, and by Theorem 
\ref{nmb:2.4}.

For $K$ compact in a Banach space $E$ and $\rho>0$, the inclusion $C^M_\rho(E\supseteq K,\mathbb{R}) \to C^L_\sigma(E\supseteq K,\mathbb{R})$ is continuous 
for all $\sigma>0$ if $M \lhd L$. It follows that the inclusion 
$\varinjlim_\rho C^M_\rho(E\supseteq K,\mathbb{R}) \to \varprojlim_\sigma C^L_\sigma(E\supseteq K,\mathbb{R})$
is continuous. This implies that the inclusion $C^{\{M\}}(U,F) \to C^{(L)}(U,F)$ is continuous 
(by definition of the structure in Section \ref{nmb:4.3}). 

Conversely, let $\mathcal{B}$ be a bounded set in $\varprojlim_{L} C^{(L)}(U,F)$, i.e., bounded in each $C^{(L)}(U,F)$. 
We claim that $\mathcal{B}$ is bounded in $C^{\{M\}}(U,F)$. We may assume without loss of generality that $E$ is a Banach space and $F=\mathbb{R}$ 
(by composing with $C^{\{M\}}(i_B,\ell)$).
Let $K\subseteq U$ be compact and $b_k:=\sup\{\|j^\infty f|_K\|_k:f\in\mathcal{B}\}$.
For all $L=(L_k)$ with $M \lhd L$ the set $\mathcal{B}$ is bounded in $C^{(L)}(U,F)$ by assumption, 
i.e., $(b_k)_k\in\bigcap_L \mathcal F^{(L)}=\mathcal F^{\{M\}}$ by Theorem \ref{nmb:2.4}.
From this follows that 
$\mathcal{B}$ is bounded in $C^{\{M\}}_K(U,F)$
and by Proposition \ref{nmb:4.2}(5) also in $C^{\{M\}}(U,F)$.

Note that this independently proves that $C^{\{M\}}(U,F)$ is $c^\infty$-complete since so is $\varprojlim_{L} C^{(L)}(U,F)$.
Moreover, it provides an independent proof of the regularity of the inductive limit involved in the definition of $C^{\{M\}}(U,F)$
if $E$ and $F$ are Banach spaces (cf.\ Proposition \ref{nmb:4.2} and the remark after Theorem \ref{nmb:8.7}). 
\end{remarks}

\begin{example}[Cartesian closedness fails without moderate growth]\label{nmb:5.4}
{\it
Let us assume that $M=(M_k)$ is weakly log-convex and has non-moderate growth 
(for instance, $M_k=q^{k^2}$, $q>1$, see \cite[2.1.3]{Thilliez08}). Then:
\begin{enumerate}
\item[\thetag{1}] 
There exists $f \in C^{\{M\}}(\mathbb{R}^2,\mathbb{C})$ such that $f^\vee :\mathbb{R} \to C^{\{M\}}(\mathbb{R},\mathbb{C})$ 
is not $C^{\{M\}}$. 
\item[\thetag{2}] 
There exists a weakly log-convex $N=(N_k)$ with $M \lhd N$ and an $f \in C^{(N)}(\mathbb{R}^2,\mathbb{C})$ such that $f^\vee :\mathbb{R} \to C^{(N)}(\mathbb{R},\mathbb{C})$ 
is not $C^{(N)}$. 
\end{enumerate}
}
\end{example}

\begin{demo}{Proof}
\thetag{1} 
There is a function $g \in C^{\{M\}}(\mathbb{R},\mathbb{C})$ such that $g^{(k)}(0) =i^k h_k$ and $h_k \ge k!\,M_k$ for all $k$; 
see \cite[Thm.\ 1]{Thilliez08}.
Defining $f(s,t):= g(s+t)$, we obtain a function $f\in C^{\{M\}}(\mathbb{R}^2,\mathbb{C})$ with 
\[
\partial^\alpha f(0,0) = i^{|\alpha|} h_{|\alpha|}, \quad h_{|\alpha|} \ge |\alpha|! \, M_{|\alpha|} \quad \text{ for all } \alpha \in \mathbb{N}^2.
\]
Since $M=(M_k)$ has non-moderate growth, there exist $j_n \nearrow \infty$ and $k_n>0$ 
such that 
\[
\Big(\frac{M_{k_n+j_n}}{M_{k_n} M_{j_n}}\Big)^\frac1{k_n+j_n} \ge n.
\]
Consider the linear functional $\ell:C^{\{M\}}(\mathbb R,\mathbb C)\to \mathbb C$
given by 
\[
\ell(g)=\sum_n\frac{i^{3 j_n}g^{(j_n)}(0)}{j_n!\,M_{j_n}\,n^{j_n}}.
\]
This functional is continuous, since 
\[
\Big|\sum_n\frac{i^{3 j_n}g^{(j_n)}(0)}{j_n!\,M_{j_n}\,n^{j_n}}\Big| 
\le 
\sum_n\frac{|g^{(j_n)}(0)|}{j_n!\, \rho^{j_n}\,M_{j_n}}\frac{\rho^{j_n}}{n^{j_n}}
\le 
C(\rho)\,\|g\|_{[-1,1], \rho} < \infty,
\]
for suitable $\rho$, where
\[
C(\rho):=\sum_n \Big(\frac{\rho}{n}\Big)^{j_n} < \infty,
\]
for all $\rho$. 
But $\ell\o f^\vee$ is not $C^{\{M\}}$, since
\begin{align*}
\|\ell\o f^\vee \|_{[-1,1],\rho_1} 
&= \sup_{k \in \mathbb N, t \in [-1,1]} 
\frac{|(\ell\o f^\vee)^{(k)}(t)|}{\rho_1^k\,k!\,M_k}
\\
&\ge
\sup_k\frac{1}{\rho_1^k\,k!\,M_k}\Big|\sum_n\frac{i^{3 j_n} f^{(j_n,k)}(0,0)}{j_n!\,M_{j_n}n^{j_n}}\Big|
\\
&= \sup_k\frac{1}{\rho_1^k\,k!\,M_k}\Big|\sum_n\frac{i^{4 j_n+k} h_{(j_n,k)}}{j_n!\,M_{j_n}n^{j_n}}\Big|
\\
&= \sup_k\frac{1}{\rho_1^k\,k!\,M_k} \sum_n\frac{h_{(j_n,k)}}{j_n!\,M_{j_n}n^{j_n}}
\\&
\ge
\sup_n\frac{1}{\rho_1^{k_n}\,k_n!\,M_{k_n}}\frac{h_{(j_n,k_n)}}{j_n!\,M_{j_n}\,n^{j_n}}
\\&
\ge
\sup_n\frac{(j_n+k_n)!\,M_{j_n+k_n}}{\rho_1^{k_n}\,k_n!\,j_n!\,M_{k_n}\,M_{j_n}\,n^{j_n}}
\ge \sup_n \frac{n^{j_n+k_n}}{\rho_1^{k_n}\,n^{j_n}} = \infty,
\end{align*}
for all $\rho_1>0$.

\thetag{2} 
By Theorem \ref{nmb:8.2}, we have, for convenient vector spaces $E,F$ and $c^\infty$-open subsets $U \subseteq E$, 
\begin{equation*} 
C^{\{M\}}(U,F) = \bigcap_N C^{(N)}(U,F),
\end{equation*}
where the intersection is taken over all weakly log-convex $N=(N_k)$ with $M \lhd N$.
Let $f$ be the function in \thetag{1}. 
Then there exist weakly log-convex sequences $N^i=(N^i_k)$, $i=1,2$, with $M \lhd N^i$ such that
$f^\vee : \mathbb{R} \to C^{(N^2)}(\mathbb{R},\mathbb{C})$ is not $C^{(N^1)}$.
By the lemma below there exists a weakly log-convex sequence $N=(N_k)$ such that $M \lhd N \le N^i$ for $i=1,2$.
Since $f \in C^{\{M\}}(\mathbb{R}^2,\mathbb{C}) \subseteq C^{(N)}(\mathbb{R}^2,\mathbb{C})$, 
the mapping $f^\vee$ has values in $C^{(N)}(\mathbb{R},\mathbb{C})$ and thus
factors over the inclusion $C^{(N)}(\mathbb{R},\mathbb{C}) \to C^{(N^2)}(\mathbb{R},\mathbb{C})$ which is obviously continuous.
It follows that $f^\vee : \mathbb{R} \to C^{(N)}(\mathbb{R},\mathbb{C})$ is not $C^{(N^1)}$ and consequently not $C^{(N)}$. 
\[
\xymatrix{
\mathbb{R} \ar@{->}[rr]^(.4){f^\vee \not\in C^{(N^1)}} \ar[rrd]_(.4){f^\vee \not\in C^{(N^1)} \supseteq C^{(N)}\quad} 
&& C^{(N^2)}(\mathbb{R},\mathbb{C}) \\
&& C^{(N)}(\mathbb{R},\mathbb{C}) \ar@{^{(}->}[u]
}
\]
By Theorem \ref{nmb:5.2}, $N=(N_k)$ has non-moderate growth.
\qed\end{demo}

\begin{lemma*}
Let $M=(M_k)$, $N^i =(N^i_k)$, $i=1,2$, be weakly log-convex with $M \lhd N^i$
for $i=1,2$. Then there exists a weakly log-convex sequence $N=(N_k)$ such that $M \lhd N \le N^i$ for $i=1,2$. 
\end{lemma*}

\begin{demo}{Proof}
Set $\bar N=(\bar N_k) := (\min\{N^1_k,N^2_k\})$ and $N=(N_k)$, 
where $(k! N_k)$ is the log-convex minorant of $(k! \bar N_k)$.
Note that $N_0=1 \le N_1$.
We have $N \le \bar N \le N^i$, 
and $M \lhd N^i$ implies $M \lhd \bar N$. 
It remains to show that $M \lhd N$.

We claim that $C^{(N)}_{\on{global}}(\mathbb{R},\mathbb{R})=C^{(\bar N)}_{\on{global}}(\mathbb{R},\mathbb{R})$, 
where for a sequence $L=(L_k) \in (\mathbb{R}_{>0})^\mathbb{N}$ we set 
\[
C^{[L]}_{\on{global}}(\mathbb{R},\mathbb{R}) := \Big\{f \in C^\infty(\mathbb{R},\mathbb{R}) : (\sup_{x\in \mathbb{R}} |f^{(k)}(x)|)_k \in \mathcal{F}^{[L]}\Big\}.
\]
In the Roumieu case this a theorem due to Cartan and Gorny, see \cite[IV E]{Koosis98};
the same proof with obvious modifications yields the Beurling version, i.e., the claim. 

Now $M \lhd \bar N$ implies 
$C^{\{M\}}_{\on{global}}(\mathbb{R},\mathbb{R}) \subseteq C^{(\bar N)}_{\on{global}}(\mathbb{R},\mathbb{R}) = C^{(N)}_{\on{global}}(\mathbb{R},\mathbb{R})$.
The function $\tilde g:=\on{Re} g+ \on{Im} g$, where $g$ is the function from the proof of \thetag{1}, 
is actually an element of $C^{\{M\}}_{\on{global}}(\mathbb{R},\mathbb{R})$ and satisfies $|\tilde g^{(k)}(0)| \ge k!\,M_k$ for all $k$; 
see \cite[Thm.\ 1]{Thilliez08}. 
Thus $\tilde g \in C^{(N)}_{\on{global}}(\mathbb{R},\mathbb{R})$ and therefore $M \lhd N$.  
\qed\end{demo}

\begin{corollary}[Canonical mappings]\label{nmb:5.5} 
Let $M=(M_k)$ be log-convex and have moderate growth.
Let $E$, $F$, etc., be convenient vector spaces and let $U$ and $V$ be
$c^\infty$-open subsets of such. 
Then we have:
\newline 
\thetag{1} The \emph{exponential law} holds: 
\[
C^{[M]}(U,C^{[M]}(V,G)) \cong C^{[M]}(U\times V, G)
\]
\indent\indent is a linear $C^{[M]}$-diffeomorphism of convenient vector spaces. 
\newline
The following canonical mappings are $C^{[M]}$.
\begin{align*}
&\operatorname{ev}: C^{[M]}(U,F)\times U\to F,\quad 
\operatorname{ev}(f,x) = f(x)
\tag{2}\\&
\operatorname{ins}: E\to C^{[M]}(F,E\times F),\quad
\operatorname{ins}(x)(y) = (x,y)
\tag{3}\\&
(\quad)^\wedge :C^{[M]}(U,C^{[M]}(V,G))\to C^{[M]}(U\times V,G)
\tag{4}\\&
(\quad)^\vee :C^{[M]}(U\times V,G)\to C^{[M]}(U,C^{[M]}(V,G))
\tag{5}\\&
\operatorname{comp}:C^{[M]}(F,G)\times C^{[M]}(U,F)\to C^{[M]}(U,G)
\tag{6}\\&
C^{[M]}(\;\;\;,\;\;\;):C^{[M]}(F,F_1)\times C^{[M]}(E_1,E)\to 
C^{[M]}\Bigl(C^{[M]}(E,F),C^{[M]}(E_1,F_1)\Bigr)
\tag{7}\\&
\qquad (f,g)\mapsto(h\mapsto f\o h\o g)
\\&
\prod:\prod C^{[M]}(E_i,F_i)\to C^{[M]}\Bigl(\prod E_i,\prod F_i\Bigr)
\tag{8}
\end{align*}
\end{corollary}

\begin{demo}{Proof}
This is a direct consequence of cartesian closedness, i.e., Theorem \ref{nmb:5.2}. See \cite[5.5]{KMRc} or 
\cite[3.13]{KM97} for the detailed arguments.
\qed\end{demo}

\section{\label{nmb:6}Uniform boundedness principles}

\begin{theorem}[$C^{[M]}$ uniform boundedness principle]\label{nmb:6.1}
Let $E$, $F$, $G$ be convenient vector spaces and let $U\subseteq F$ be
$c^\infty$-open. A linear mapping $T:E\to C^{[M]}(U,G)$
is bounded if and only if $\on{ev}_x\o T: E\to G$ is bounded for every $x\in U$. 
\end{theorem}

\begin{demo}{Proof} 
$(\Rightarrow)$ 
For $x\in U$ and $\ell\in G^*$, the linear mapping $\ell\o\on{ev}_x = C^{[M]}(x,\ell):C^{[M]}(U,G)\to
\mathbb{R}$ is continuous, thus $\on{ev}_x$ is bounded. 
Therefore, if $T$ is bounded then so is $\on{ev}_x\o T$.

$(\Leftarrow)$
Suppose that $\on{ev}_x\o T$ is bounded for all $x\in U$.
By the definition of $C^{[M]}(U,G)$ in Section \ref{nmb:4.3}
it is enough to show that $T$ is bounded in the case that $E$ and $F$ are Banach spaces  and $G=\mathbb R$.
By Section \ref{nmb:4.1}, $C^{[M]}(U,\mathbb{R})=\varprojlim_{K} C^{[M]}(F\supseteq K,\mathbb{R})$,
by Proposition \ref{nmb:4.2}(2), $C^{(M)}(F\supseteq K,\mathbb{R})$ is a Fr\'echet space, and by Proposition \ref{nmb:4.2}(3),
$C^{\{M\}}(F\supseteq K,\mathbb{R})$ is an (LB)-space, so $C^{[M]}(F\supseteq K,\mathbb{R})$ is webbed and hence the 
closed graph theorem \cite[52.10]{KM97} gives the desired result.
\begin{equation}\label{eq:27}
\begin{split}
\xymatrix{
E\ar@{->}[r]^(0.4){T} \ar@{->}[dr] &C^{[M]}(U,\mathbb{R}) \ar@{=}[d] 
\ar@{->}[r]^(0.6){\on{ev}_x} &\mathbb{R} \\
&\varprojlim _{K}C^{[M]}(F\supseteq K,\mathbb{R})\ar@{->}[r] &
C^{[M]}(F\supseteq K,\mathbb{R} )\ar@{->}[u]_{(\on{ev}_x)^0} 
}
\end{split}
\end{equation}
\qed\end{demo}

\begin{remark}\label{nmb:6.2}
Alternatively, the $C^{\{M\}}$ uniform boundedness principle follows from the $C^{(M)}$ uniform boundedness principle 
and from the remark after Theorem \ref{nmb:8.2}, since the structure of $C^{\{M\}}(U,F)=\varprojlim_{L} C^{(L)}(U,F)$ 
is initial with respect to the inclusions $\varprojlim_{L} C^{(L)}(U,F) \to C^{(L)}(U,F)$ for all $L$.
This is no circular argument, since the first identity in Theorem \ref{nmb:8.2} was proved in Remark \ref{nmb:5.3} without using the 
$C^{\{M\}}$ uniform boundedness principle, i.e., Theorem \ref{nmb:6.1}.
\end{remark}

\section{\label{nmb:7}Relation to previously considered structures}

In \cite{KMRc} and \cite{KMRq} we have developed the convenient setting for all \emph{reasonable} non-quasianalytic and some quasianalytic 
(namely, $\mathcal{L}$-intersectable, see Section \ref{nmb:7.1}) 
Denjoy--Carleman classes of Roumieu type. 
We have worked with a definition which is based on \emph{testing along curves}. The resulting structures were denoted by $C^M$ in \cite{KMRc} and \cite{KMRq} and will 
be denoted by $C_{\on{curve}}^{\{M\}}$ in this section; this notation does not appear elsewhere in this paper.
We shall now show that they coincide bornologically with the structure $C^{\{M\}}$ studied in the present paper.
Furthermore, we prove that the bornologies induced by $C^{\{1\}}$ and the structure $C^\omega$ of real analytic mappings introduced in \cite{KrieglMichor90} are isomorphic; here 
$1$ denotes the constant sequence $(1)_k$. Note that $C^{\{1\}}$ is not $\mathcal{L}$-intersectable (see 
\cite[1.8]{KMRq}).

\subsection{\label{nmb:7.1}Testing along curves}
Let $M=(M_k)$ be log-convex, $E$ and $F$ convenient vector spaces, and $U$ a $c^\infty$-open subset in $E$.
If $M=(M_k)$ is non-quasianalytic we set
\[
C_{\on{curve}}^{\{M\}}(U,F) := \Bigl\{f \in F^U:
\forall \ell \in F^*\; \forall c \in C^{\{M\}}(\mathbb{R},U): \ell \o f \o c \in C^{\{M\}}(\mathbb{R},\mathbb{R})\Bigr\}.
\]
If $M=(M_k)$ is quasianalytic and \emph{$\mathcal{L}$-intersectable}, i.e., $\mathcal{F}^{\{M\}} = \bigcap_{L \in \mathcal{L}(M)} \mathcal{F}^{\{L\}}$, 
where 
\[
\mathcal{L}(M) := \Bigl\{L=(L_k) : L \ge M, L \text{ is non-quasianalytic log-convex}\Bigr\},
\]
we define 
\[
C_{\on{curve}}^{\{M\}}(U,F) := \bigcap_{L \in \mathcal{L}(M)} C_{\on{curve}}^{\{L\}}(U,F).
\]
Note that non-quasianalytic log-convex sequences are trivially $\mathcal{L}$-intersectable. For 
non-quasianalytic $M=(M_k)$
we supply $C_{\on{curve}}^{\{M\}}(U,F)$ with the initial locally convex structure induced by all linear mappings:
\begin{align*}
C_{\on{curve}}^{\{M\}}(U,F) &\East{C_{\on{curve}}^{\{M\}}(c,\ell)}{} C^{\{M\}}(\mathbb{R}, \mathbb R), \quad f \mapsto \ell\o f\o c, \qquad \ell \in F^*, c \in C^{\{M\}}(\mathbb{R},U),
\end{align*}
and for quasianalytic and \emph{$\mathcal{L}$-intersectable} $M=(M_k)$
by all inclusions 
\begin{align*}
C_{\on{curve}}^{\{M\}}(U,F) &\longrightarrow C_{\on{curve}}^{\{L\}}(U, F), \qquad L \in \mathcal{L}(M). 
\end{align*}
In both cases $C_{\on{curve}}^{\{M\}}(U,F)$ is a convenient vector space.

Let $C^\omega(\mathbb{R},\mathbb{R})$ denote the real analytic functions $f : \mathbb{R} \to \mathbb{R}$ and set
\begin{align*}
C^\omega(U,F) &:= \Bigl\{f \in C^\infty(U,F) : \forall \ell \in F^*\; \forall c \in C^\omega(\mathbb{R},U): 
\ell \o f \o c \in C^\omega(\mathbb{R},\mathbb{R})\Bigr\},
\end{align*}
where $C^\omega(\mathbb R,U)$ is the space of all weakly $C^\omega$-curves in $U$.
We equip $C^\omega(U,\mathbb{R})$ with the initial locally convex structure induced by the family of mappings
\begin{align*} 
C^\omega(U,\mathbb{R}) &\East{{c^*}}{} C^\omega(\mathbb R,\mathbb R), \quad f \mapsto f\o c,\qquad c\in C^\omega(\mathbb{R},U),\\ 
C^\omega(U,\mathbb{R}) &\East{{c^*}}{} C^\infty(\mathbb R,\mathbb R), \quad f \mapsto f\o c,\qquad c\in C^\infty(\mathbb R,U),
\end{align*}
where $C^\infty(\mathbb R,\mathbb R)$ carries the topology of compact 
convergence in each derivative separately, and 
where $C^\omega(\mathbb R,\mathbb R)$ is equipped with the final locally 
convex topology 
with respect to the embeddings (restriction mappings) of all spaces 
of holomorphic mappings from a neighborhood $V$ of $\mathbb R$ in 
$\mathbb C$ mapping $\mathbb R$ to $\mathbb R$, and each of these spaces 
carries the topology of compact convergence.
The space
$C^\omega(U,F)$ is equipped with the initial locally convex structure induced by all mappings 
\[
C^\omega(U,F) \East{{\ell_*}}{}C^\omega(U,\mathbb R), \quad f \mapsto \ell \o f, \qquad \ell \in F^*.
\]
This is again a convenient vector space.

\begin{theorem} \label{nmb:7.2}
Let $M=(M_k)$ be log-convex, $E$ and $F$ convenient vector spaces, and $U$ a $c^\infty$-open subset in $E$.
We have:
\begin{enumerate}
\item[\thetag{1}] If $M=(M_k)$ is $\mathcal{L}$-intersectable, then
\[
C^{\{M\}}(U,F) = C_{\on{curve}}^{\{M\}}(U,F)
\] 
as vector spaces with bornology. 
\item[\thetag{2}] If $1$ denotes the constant sequence, then 
\[
C^{\{1\}}(U,F) = C^\omega(U,F)
\] 
as vector spaces with bornology.
\end{enumerate}
\end{theorem}

\begin{demo}{Proof}
\thetag{1} If $M=(M_k)$ is non-quasianalytic, then $C^{\{M\}}(U,F)$ and $C_{\on{curve}}^{\{M\}}(U,F)$ coincide as vector spaces, by \cite[2.8]{KMRq}.
If $M=(M_k)$ is quasianalytic and $\mathcal{L}$-intersectable, then the non-quasianalytic case implies that
\[
C_{\on{curve}}^{\{M\}}(U,F) = \bigcap_{L \in \mathcal{L}(M)} C_{\on{curve}}^{\{L\}}(U,F) = \bigcap_{L \in \mathcal{L}(M)} C^{\{L\}}(U,F) = C^{\{M\}}(U,F)
\]
as vector spaces, where the last equality is a consequence of the definition of $C^{\{M\}}(U,F)$ (see Section \ref{nmb:4.3}) and 
of \cite[1.6]{KMRq} (applied to $C^{\{M\}}(U_B,\mathbb{R})$). 
The fact that both spaces $C^{\{M\}}(U,F)$ and $C_{\on{curve}}^{\{M\}}(U,F)$ are convenient and 
satisfy the uniform boundedness principle with respect to the set of point evaluations, 
see Theorem \ref{nmb:6.1} and \cite[2.9]{KMRq}, 
implies that the identity is a bornological isomorphism.

\thetag{2} We show first that $C^{\{1\}}(U,F) = C^\omega(U,F)$ as vector spaces. By Section \ref{nmb:4.3} and \cite[10.6]{KM97}, 
it suffices to consider the case that $U$ is open in a Banach space $E$ and $F=\mathbb{R}$.

Let $f \in C^\omega(U,\mathbb{R})$. By \cite[2.4 and 2.7]{KrieglMichor90} or \cite[10.1 and 10.4]{KM97}, 
this is equivalent to $f$ being smooth and being locally given by its convergent Taylor series.
Let $K \subseteq U$ be compact. Since the Taylor series of $f$ converges locally, there exist constants $C,\rho>0$ such that 
\begin{equation*} 
\frac{\|f^{(k)}(a)\|_{L^k(E,\mathbb{R})}}{k!} \le C \rho^k, \quad \text{ for all } a \in K, k \in \mathbb{N},
\end{equation*} 
that is, $f \in C^{\{1\}}(U,\mathbb{R})$.
\newline
Conversely, the above estimate for compact subsets $K$ of affine lines in $E$  
implies that the restriction of $f$ to each affine line is real analytic and hence 
$f\in C^\omega(U,\mathbb R)$ by \cite[10.1]{KM97}.

The bornologies coincide, since both spaces are convenient and satisfy the uniform boundedness 
principle with respect to the set of point evaluations, see Theorem \ref{nmb:6.1} and 
\cite[5.6]{KrieglMichor90} or \cite[11.12]{KM97}.
\qed\end{demo}

\section{\label{nmb:8}More on function spaces}

\begin{proposition}[Inclusions]\label{nmb:8.1} 
Let $M=(M_k)$, $N=(N_k)$ be positive sequences, 
$E$, $F$ convenient vector spaces, and $U \subseteq E$ a $c^\infty$-open subset. We have:
\begin{enumerate}
\item[\thetag{1}] 
$C^{(M)}(U,F) \subseteq C^{\{M\}}(U,F)\subseteq C^\infty(U,F)$. 
\item[\thetag{2}] If there exist $C,\rho>0$ so that $M_k \le C \rho^k N_k$ for all $k$, then 
\[
C^{(M)}(U,F)\subseteq C^{(N)}(U,F) \quad \text{ and } \quad C^{\{M\}}(U,F)\subseteq C^{\{N\}}(U,F).
\]
\item[\thetag{3}] If for each $\rho>0$ there exists $C>0$ so that $M_k \le C \rho^k N_k$ for all $k$, i.e., $M \lhd N$, then
\[
C^{\{M\}}(U,F) \subseteq C^{(N)}(U,F).
\]
\item[\thetag{4}] For $U \ne \emptyset$ and $F\ne\{0\}$ we have: 
\begin{align*}
C^\omega(U,F) \subseteq C^{(M)}(U,F) &\Longleftrightarrow M_k^{1/k} \to \infty, \quad \text{and}\\
C^\omega(U,F) \subseteq C^{\{M\}}(U,F) &\Longleftrightarrow \varliminf M_k^{1/k} > 0. 
\end{align*}
\end{enumerate}
All these inclusions are bounded.
\end{proposition}

\begin{demo}{Proof}
The inclusions in \thetag{1}, \thetag{2}, and \thetag{3} follow immediately from the definitions in Sections \ref{nmb:4.1} and 
\ref{nmb:4.3} and Lemma \ref{nmb:2.3}. Here we use that $C^{\{1\}}(U,F) = C^\omega(U,F)$ as
vector spaces with bornology,
see Theorem \ref{nmb:7.2}.

The directions ($\Leftarrow$) in \thetag{4} are direct consequences of \thetag{2} and \thetag{3}.
The directions ($\Rightarrow$) follow, since they have been shown in Section \ref{nmb:2.1} for 
$E=F=\mathbb R$.

All inclusions are bounded, since all spaces are convenient and satisfy the uniform boundedness principle, cf.\  
Theorem \ref{nmb:6.1} and \cite[5.26]{KM97} for $C^\infty$.
\qed\end{demo}

\begin{theorem} \label{nmb:8.2}
Let $M=(M_k)$ be (weakly) log-convex, $E$ and $F$ convenient vector spaces, and $U$ a $c^\infty$-open subset in $E$. 
We have 
\[
C^{\{M\}}(U,F) = \varprojlim_{L} C^{(L)}(U,F) = \varprojlim_{L} C^{\{L\}}(U,F)
\]
as vector spaces with bornology,
where the projective limits are taken over all (weakly) log-convex sequences $L=(L_k)$ with $M \lhd L$.
\end{theorem}

\begin{demo}{Proof}
The three spaces coincide as vector spaces: By Section \ref{nmb:4.3} it suffices to 
assume that $E$ and $F$ are Banach spaces, and by Section \ref{nmb:4.1} and Proposition \ref{nmb:4.2}(5) it suffices to apply Theorem 
\ref{nmb:2.4} to the sequence
$(\|j^\infty f|_K\|_m)$.

Each space is convenient (see Section \ref{nmb:4.3}; projective limits preserve $c^\infty$-completeness)
and each space satisfies the uniform boundedness principle with respect to the set of point evaluations
(see Theorem \ref{nmb:6.1}; the structure of $\varprojlim_{L} C^{[L]}(U,F)$ 
is initial with respect to the inclusions $\varprojlim_{L} C^{[L]}(U,F) \to C^{[L]}(U,F)$ for all $L$).
Thus the identity between any two of the three spaces is a bornological isomorphism. 
\qed\end{demo}

\begin{remark*}
By the remark after Theorem \ref{nmb:2.4} the statement of the theorem still holds, if $M=(M_k)$ is just any positive sequence,
where the projective limits are now taken over all positive sequences $L=(L_k)$ with $M \lhd L$. 
\end{remark*}

\begin{proposition}[Derivatives]\label{nmb:8.3} 
Let $M=(M_k)$ be a positive sequence and set $M_{+1} = (M_{k+1})$. 
Let $E$ and $F$ be convenient vector spaces, and $U\subseteq E$ a $c^\infty$-open subset.
Then we have:
\begin{enumerate}
\item[\thetag{1}] Multilinear mappings between convenient vector spaces are $C^{[M]}$ if and only
if they are bounded.
\item[\thetag{2}] If $f:E\supseteq U\to F$ is $C^{[M]}$, then the derivative 
$df:U\to L(E,F)$ is $C^{{[M_{+1}]}}$, 
where the space $L(E,F)$ of all bounded linear mappings is considered 
with the topology of uniform convergence on bounded sets.
If $M_{+1}=(M_{k+1})$ is weakly log-convex (which is the case if $M=(M_k)$ is weakly log-convex), also $(df)^\wedge :U\times E\to F$ is 
$C^{{[M_{+1}]}}$,
\item[\thetag{3}] The chain rule holds.
\end{enumerate}
\end{proposition}

\begin{demo}{Proof}
\thetag{1}
If $f$ is $C^{[M]}$ then it is smooth and hence
bounded by \cite[5.5]{KM97}. 
Conversely, if $f$ is multilinear and bounded then
it is smooth, again by \cite[5.5]{KM97}. 
Further\-more, $f\o i_B$ is multilinear and continuous and all derivatives of high order
vanish. Thus $f$ is $C^{[M]}$, by Section \ref{nmb:4.3}. 

\thetag{2} Since $f$ is smooth, by \cite[3.18]{KM97} 
the mapping $df:U\to L(E,F)$ exists and is smooth. 
We have to show that $(df) \o i_B : U_B \to L(E,F)$ is $C^{{[M_{+1}]}}$, 
for all closed absolutely convex bounded subsets $B \subseteq E$. 
By the uniform boundedness principle 
\cite[5.18]{KM97} and by Lemma \ref{nmb:5.1} it suffices to show that the mapping 
$U_B \ni x \mapsto \ell(df(i_B(x))(v))\in \mathbb{R}$ is $C^{{[M_{+1}]}}$ for each
$\ell\in F^*$ and $v\in E$.

Since $\ell\o f$ is $C^{(M)}$ (resp.\ $C^{\{M\}}$), for each closed absolutely convex bounded $B \subseteq E$, 
each compact $K \subseteq U_B$, and each $\rho>0$
(resp. some $\rho>0$)
the set
\[
\Big\{\frac{\|d^k(\ell\o f\o i_B)(a)\|_{L^k(E_B,\mathbb{R})}}{k!\,\rho^k\, M_k} : a \in K,k \in \mathbb{N}\Big\}
\]
is bounded, say by $C>0$. The assertion follows in both cases from the following computation. 
For $v\in E$ and those $B$ containing $v$ we then have:
\begin{align*}
\|d^k(L(\ell,v)\o &df)\o i_B)(a)\|_{L^k(E_B,\mathbb R)} 
=\|d^k(d(\ell\o f)(\quad)(v))\o i_B)(a)\|_{L^k(E_B,\mathbb R)}\\
&=\|d^{k+1}(\ell\o f\o i_B)(a)(v,\dots)\|_{L^k(E_B,\mathbb R)} 
\\&
\le\|d^{k+1}(\ell\o f\o i_B)(a)\|_{L^{k+1}(E_B,\mathbb R)}\|v\|_{B} 
\le C\, (k+1)!\, \rho^{k+1}\, M_{k+1}
\\&
= C\rho\,((k+1)^{1/k}\rho)^{k}\, k!\, M_{k+1}
\leq C\rho\,(2\rho)^k\,k!\,(M_{+1})_k.
\end{align*}
By Proposition \ref{nmb:8.4} below also $(df)^\wedge$ is $C^{{[M_{+1}]}}$, if $M=(M_k)$ is weakly log-convex.

\thetag{3} This is valid even for all smooth $f$ by \cite[3.18]{KM97}. 
\qed\end{demo}

\begin{proposition} \label{nmb:8.4} 
We have:
\begin{enumerate}
\item[\thetag{1}] For convenient vector spaces $E$ and $F$, the following topologies have the same bounded subsets
in $L(E,F)$:
\begin{itemize}
\item The topology of uniform convergence on bounded subsets
of $E$.
\item The topology of pointwise convergence.
\item The trace topology of $C^\infty(E,F)$.
\item The trace topology of $C^{[M]}(E,F)$.
\end{itemize}
\item[\thetag{2}] Let $M=(M_k)$ be weakly log-convex, $E$, $F$, and $G$ convenient vector spaces, and $U\subseteq E$ a 
$c^\infty$-open subset. A mapping $f:U\times F\to G$ which is linear in
the second variable is $C^{[M]}$ if and only if 
$f^\vee:U\to L(F,G)$ is well defined and $C^{[M]}$.
\end{enumerate}
Analogous results hold for spaces of multilinear mappings. 
\end{proposition}

\begin{demo}{Proof}
\thetag{1}
That the first three topologies on $L(E,F)$ have the same bounded sets has been shown
in \cite[5.3 and 5.18]{KM97}.
The inclusion $C^{[M]}(E,F)\to C^\infty(E,F)$ is bounded by Proposition \ref{nmb:8.1}.
Conversely, the inclusion $L(E,F)\to C^{[M]}(E,F)$ is bounded
by the uniform boundedness principle, i.e., Theorem \ref{nmb:6.1}.

\thetag{2} The assertion for $C^\infty$ is true by \cite[3.12]{KM97} 
since $L(E,F)$ is closed in $C^\infty(E,F)$.

Suppose that $f$ is $C^{[M]}$. We have to show 
that $f^\vee \o i_B$ is $C^{[M]}$ into $L(F,G)$, for all closed absolutely convex bounded subsets $B \subseteq E$. By the uniform boundedness principle 
\cite[5.18]{KM97} and by Lemma \ref{nmb:5.1} it suffices to show that the mapping
$U_B \ni x \mapsto \ell\bigl(f^\vee(i_B(x))(v)\bigr)=\ell\bigl(f(i_B(x),v)\bigr)\in \mathbb R$ is $C^{[M]}$ for each
$\ell\in G^*$ and $v\in F$; this is obviously true. 

Conversely, let $f^\vee:U\to L(F,G)$ be $C^{[M]}$. By \thetag{1} the inclusion $L(F,G) \to C^{[M]}(F,G)$ is bounded linear, 
and so $f^\vee:U\to C^{[M]}(F,G)$ is $C^{[M]}$.
By cartesian closedness, i.e., Theorem \ref{nmb:5.2} (the direction which holds without moderate growth), 
$f:U\times F\to G$ is $C^{[M]}$ and linearity in the second variable is obvious. 
\qed\end{demo}

\begin{remark*}
We may prove 
$f^\vee \in C^{[M]}(U,L(F,G)) \Rightarrow f \in C^{[M]}(U \times F,G)$ without using cartesian closedness: 
By composing with $\ell\in G^*$ we
may assume that $G=\mathbb R$.
By induction we have: 
\begin{multline*}
d^kf(x,w_0)\big((v_k,w_k),\dots,(v_1,w_1)\big) 
= d^k (f^\vee)(x)(v_k,\dots,v_1)(w_0) +
\\
+ \sum_{i=1}^k d^{k-1}(f^\vee)(x)(v_k,\dots,\widehat{v_i,}\dots,v_1)(w_i)
\end{multline*}
Thus for $B$, $B'$ closed absolutely convex bounded subsets of $E$, $F$, respectively, $K \subseteq U_B$ compact, and $x \in K$ we have: 
\begin{align*}
&\|d^k f (x,w_0)\|_{L^k(E_B\times F_{B'},\mathbb R)}\le
\\&
\le \|d^k (f^\vee)(x)(\dots)(w_0)\|_{L^k(E_B,\mathbb R)}
+ \sum_{i=1}^k \|d^{k-1}(f^\vee)(x)\|_{L^{k-1}(E_B,L(F_{B'},\mathbb R))}
\\&
\le \|d^k (f^\vee)(x)\|_{L^k(E_B,L(F_{B'},\mathbb R))}\|w_0\|_{B'}
+ k\, \|d^{k-1}(f^\vee)(x)\|_{L^{k-1}(E_B,L(F_{B'},\mathbb R))}
\\&
\le C\, \rho^k\, k!\, M_k \|w_0\|_{B'}
+ k\, C \,\rho^{k-1}\, (k-1)!\, M_{k-1} 
= C\,\rho^k \,k!\,M_k \Bigl(\|w_0\|_{B'}+ \tfrac{M_{k-1}}{\rho\, M_k}\Bigr),
\end{align*}
for all $\rho>0$ and some $C=C(\rho)$ (resp.\ for some $C,\rho>0$), 
since the mapping $L(i_{B'},\mathbb R)\o f^\vee \o i_B: U_B\to L(F_{B'},\mathbb R)$ is $C^{[M]}$.
Since $k \mapsto k!\, M_k$ is increasing (see the remarks in Section \ref{nmb:2.6}), we have $\frac{M_{k-1}}{M_k} \le k \le 2^k$, 
and we may conclude that 
$f$ is $C^{[M]}$.
\end{remark*}

Let $r=(r_k)$ be a positive sequence,
$E$ and $F$ Banach spaces,
and $K\subseteq E$ compact convex. 
Consider 
\begin{align*}
C^{M}_{(r_k)}(E\supseteq K,F) :&= 
\Bigl\{(f^m)_m\in \prod_{m\in\mathbb N}C(K,L_{\text{sym}}^m(E,F)):\|f\|_{(r_k)}<\infty\Bigr\},\\
\intertext{where }\|f\|_{(r_k)}:&=\max\biggl\{\sup\Bigl\{\frac{\|f\|_m}{m!\, r_m\, M_m}:m\in\mathbb N\Bigr\},\\
&\qquad\qquad\sup\Bigl\{\frac{|||f|||_{n,k}}{(n+k+1)!\,r_{n+k+1}\,M_{n+k+1}}:k,n\in\mathbb N\Bigr\}\biggr\}.
\end{align*}
If $(r_k)=(\rho^k)$ for some $\rho>0$ we just write $\rho$ instead of $(r_k)$ as indices and recover the spaces introduced in Section \ref{nmb:4.1}.
Similarly as in Proposition \ref{nmb:4.2}(1) one shows that
the spaces $C^M_{(r_k)}(E \supseteq K,F)$ are Banach spaces.

\begin{theorem} \label{nmb:8.6} 
Let $E$ and $F$ be Banach spaces and let $U\subseteq E$ be open and convex.
Then we have 
\[
C^{(M)}(U,F) = \varprojlim_{K,(r_k)} C^M_{(r_k)}(E\supseteq K,F)
\]
as vector spaces with bornology.
Here $K$ runs through all compact convex subsets of $U$ ordered by inclusion 
and
$(r_k)$ runs through all sequences of positive real numbers for which $\rho^k/r_k\to 0$ for some $\rho>0$.
\end{theorem}

\begin{demo}{Proof} 
Note first that the elements of the space $\varprojlim_{K,(r_k)} C^M_{(r_k)}(E\supseteq K,F)$ are
smooth functions $f:U\to F$ which can be seen as in the proof of Proposition \ref{nmb:4.2}(5).
By Lemma \ref{nmb:4.7} it coincides with $C^{(M)}(U,F)$ as vector space.

Obviously the identity is continuous from left to right.
The space on the right-hand side is as a projective limit of Banach spaces convenient and $C^{(M)}(U,F)$
satisfies the uniform boundedness principle, i.e., Theorem \ref{nmb:6.1}, with respect to the set of point evaluations.
Thus the identity from right to left is bounded.
\qed\end{demo}

\begin{theorem} \label{nmb:8.7} 
Let $E$ and $F$ be Banach spaces and let $U\subseteq E$ be open and convex.
Then we have 
\[
C^{\{M\}}(U,F) = \varprojlim_{K,(r_k)} C^M_{(r_k)}(E\supseteq K,F)
\]
as vector spaces with bornology.
Here $K$ runs through all compact convex subsets of $U$ ordered by inclusion 
and
$(r_k)$ runs through all sequences of positive real numbers for which $\rho^k/r_k\to 0$ for all 
$\rho>0$.
\end{theorem}

\begin{demo}{Proof} 
The proof is literally identical with the proof of Theorem \ref{nmb:8.6}, where we replace $C^{(M)}$ with $C^{\{M\}}$ and use 
Lemma \ref{nmb:4.8} instead of Lemma \ref{nmb:4.7}.
\qed\end{demo}

\begin{remark*}
Let us prove that the identity $\varprojlim_{K,(r_k)} C^M_{(r_k)}(E\supseteq K,F) \to C^{\{M\}}(U,F)$ is bounded without using 
the $C^{\{M\}}$ uniform boundedness principle, i.e., Theorem \ref{nmb:6.1}:
Let $\mathcal{B}$ be a bounded set in $\varprojlim_{K,(r_k)} C^M_{(r_k)}(E\supseteq K,F)$, i.e., for each compact $K$ and each $(r_k)$ with $\rho^k/r_k\to 0$ for all 
$\rho>0$ the set $\mathcal{B}$ is bounded in $C^M_{(r_k)}(E \supseteq K,F)$, i.e.,
\[
\sup\{\|f|_K\|_{(r_k)} : f \in \mathcal{B}\} < \infty.
\] 
Since the elements of $\varprojlim_{K,(r_k)} C^M_{(r_k)}(E\supseteq K,F)$ are the infinite jets of smooth functions, we may estimate $|||f|_K|||_{n,k}$ by $\|f|_K\|_{n+k+1}$ by \eqref{eq:12},
and so the sequence
\[
a_k := \sup\Big\{\frac{\|f|_K\|_k}{k!\, M_k} : f \in \mathcal{B}\Big\} < \infty
\] 
satisfies $\sup_{k} a_k/r_k < \infty$ for each $(r_k)$ as above. 
By \cite[9.2]{KM97}, these are the coefficients of a power series with positive radius of convergence.
Thus $a_k/\rho^k$ is bounded for some $\rho>0$. That means that $\mathcal{B}$ is contained and bounded in $C^M_{\rho}(E \supseteq K,F)$.

This also provides an independent proof of the completeness of $C^{\{M\}}(U,F)$ and of the regularity of the involved 
inductive limit (cf.\ Proposition \ref{nmb:4.2}
and Remark \ref{nmb:5.3}).
\end{remark*}

\begin{lemma} \label{nmb:8.8}
For convenient
vector spaces $E$, $F$, $G$, and $c^\infty$-open $V\subseteq F$  
the flip of variables induces an isomorphism $L(E,C^{[M]}(V,G)) \cong C^{[M]}(V,L(E,G))$ as vector spaces. 
\end{lemma}

\begin{demo}{Proof} 
For $f \in C^{[M]}(V,L(E,G))$ consider $\tilde f(x) := \on{ev}_x\o f\in C^{[M]}(V,G)$ for $x\in E$. 
By the uniform boundedness principle, i.e., Theorem \ref{nmb:6.1}, 
the linear mapping $\tilde f$ is bounded, since 
$\on{ev}_y\o\tilde f=f(y)\in L(E,G)$ for $y \in V$.

If conversely $\ell\in L(E,C^{[M]}(V,G))$, we consider
$\tilde \ell(y) = \on{ev}_y\o \ell \in L(E,G)$ for $y\in V$. 
Since the bornology of $L(E,G)$ (see Proposition \ref{nmb:8.4}) is generated by $\mathcal{S}:=\{\on{ev}_x:x\in E\}$ and since 
$\on{ev}_x\o\tilde \ell=\ell(x)\in C^{[M]}(V,G)$, it follows that 
$\tilde \ell:V \to L(E,G)$ is $C^{[M]}$, by Lemma \ref{nmb:5.1} (and by composing with all $i_B : V_B \to V$).
\qed\end{demo}

\begin{lemma} \label{nmb:8.9}
Let $E$ be a convenient vector space and let $U \subseteq E$ be $c^\infty$-open.
By $\lambda^{[M]}(U)$ we denote the $c^\infty$-closure of the linear subspace
generated by $\{\on{ev}_x: x\in U\}$ in $C^{[M]}(U,\mathbb R)'$ and let $\delta:U\to\lambda^{[M]}(U)$
be given by $x\mapsto \on{ev}_x$.
Then $\lambda^{[M]}(U)$ is the free convenient vector space over $C^{[M]}$, i.e.,
for every convenient vector space $G$ the $C^{[M]}$-mapping $\delta$ induces a
bornological isomorphism
\begin{align*}
L(\lambda^{[M]}(U),G)&\cong C^{[M]}(U,G).
\end{align*}
\end{lemma}

\begin{demo}{Proof}
The proof goes along the same lines as in \cite[5.1.1]{FK88} and \cite[23.6]{KM97}.
Note first that $\lambda^{[M]}(U)$ is a convenient vector space, since it is $c^\infty$-closed in the convenient vector space 
$C^{[M]}(U,\mathbb R)'$. Moreover, $\delta$ is $C^{[M]}$, by Lemma \ref{nmb:5.1} (and by composing with all $i_B : U_B \to U$), 
since $\on{ev}_h\o\delta=h$ for all $h\in C^{[M]}(U,\mathbb R)$. So $\delta^*:L(\lambda^{[M]}(U),G)\to C^{[M]}(U,G)$
is a well-defined linear mapping. 
This mapping is injective, since each bounded linear mapping $\lambda^{[M]}(U)\to G$ is uniquely 
determined on $\delta(U)=\{\on{ev}_x: x\in U\}$.
Let now $f\in C^{[M]}(U,G)$. Then 
$\ell\o f\in C^{[M]}(U,\mathbb R)$ for every $\ell\in G^*$, and hence 
$\tilde f:C^{[M]}(U,\mathbb R)'\to \prod_{G^*}\mathbb R$ given by $\tilde f(\phi)=(\phi(\ell\o f))_{\ell\in G^*}$
is a well-defined bounded linear mapping. Since it maps $\on{ev}_x$ to $\tilde f(\on{ev}_x)=\delta(f(x))$, 
where $\delta:G\to\prod_{G^*}\mathbb R$ denotes the bornological embedding given by 
$y\mapsto (\ell(y))_{\ell\in G^*}$, it induces a bounded linear mapping $\tilde f:\lambda^{[M]}(U)\to G$
satisfying $\tilde f\o\delta=f$.
Thus $\delta^*$ is a linear bijection.
That it is a bornological isomorphism follows from 
the uniform boundedness principle, i.e., Theorem \ref{nmb:6.1}, and from Proposition \ref{nmb:8.4}.
\qed\end{demo}

\begin{theorem}[Canonical isomorphisms]\label{nmb:8.10}
Let $M=(M_k)$ and $N=(N_k)$ be positive sequences.
Let $E$, $F$ be convenient vector spaces and let
$W_i$ be $c^\infty$-open subsets in such. 
We have the following natural bornological isomorphisms:
\begin{enumerate}
\item[\thetag{1}] $C^{(M)}(W_1,C^{(N)}(W_2,F))\cong C^{(N)}(W_2,C^{(M)}(W_1,F))$,
\item[\thetag{2}] $C^{\{M\}}(W_1,C^{\{N\}}(W_2,F))\cong C^{\{N\}}(W_2,C^{\{M\}}(W_1,F))$,
\item[\thetag{3}] $C^{(M)}(W_1,C^{\{N\}}(W_2,F))\cong C^{\{N\}}(W_2,C^{(M)}(W_1,F))$,
\item[\thetag{4}] $C^{[M]}(W_1,C^\infty(W_2,F))\cong C^\infty(W_2,C^{[M]}(W_1,F))$.
\item[\thetag{5}] $C^{[M]}(W_1,C^\omega(W_2,F))\cong C^\omega(W_2,C^{[M]}(W_1,F))$.
\item[\thetag{6}] $C^{[M]}(W_1,L(E,F))\cong L(E,C^{[M]}(W_1,F))$.
\item[\thetag{7}] $C^{[M]}(W_1,\ell^\infty(X,F))\cong \ell^\infty(X,C^{[M]}(W_1,F))$.
\item[\thetag{8}] $C^{[M]}(W_1,\operatorname{\mathcal Lip}^k(X,F))\cong \operatorname{\mathcal Lip}^k(X,C^{[M]}(W_1,F))$.
\end{enumerate}
In \thetag{7} the space $X$ is an $\ell^\infty$-space, i.e., 
a set together 
with a bornology induced by a family of real valued functions on $X$, 
cf.\ \cite[1.2.4]{FK88}.
In \thetag{8} the space $X$ is a $\operatorname{\mathcal Lip}^k$-space, cf.\ 
\cite[1.4.1]{FK88}.
The spaces $\ell^\infty(X,F)$ and $\operatorname{\mathcal Lip}^k(X,F)$ are defined in 
\cite[3.6.1 and 4.4.1]{FK88}. 
\end{theorem}

\begin{demo}{Proof} 
Let $\mathcal C^1$ and $\mathcal C^2$ denote any of the functions spaces mentioned above
and $X_1$ and $X_2$ the corresponding domains.
In order to show that the flip of coordinates $f\mapsto \tilde f$, 
$\mathcal C^1(X_1,\mathcal C^2(X_2,F))\to \mathcal C^2(X_2,\mathcal C^1(X_1,F))$
is a well-defined bounded linear mapping we have to show:
\begin{itemize}
\item $\tilde f(x_2)\in\mathcal C^1(X_1,F)$, which is obvious, since
$\tilde f(x_2)=\on{ev}_{x_2}\o f:X_1\to \mathcal C^2(X_2,F)\to F$.
\item $\tilde f\in\mathcal C^2(X_2,\mathcal C^1(X_1,F))$, which we will show below.
\item $f\mapsto \tilde f$ is bounded and linear, which follows by applying the appropriate
uniform boundedness theorems for $\mathcal C^2$ and $\mathcal C^1$, since
$f\mapsto \on{ev}_{x_1}\o \on{ev}_{x_2}\o\tilde f=\on{ev}_{x_2}\o\on{ev}_{x_1}\o f$ is bounded and linear.
\end{itemize}
All occurring 
function spaces are convenient and satisfy the uniform 
$\mathcal{S}$-boundedness theorem, where $\mathcal{S}$ is the set of point 
evaluations:
\begin{itemize}
\item[$C^{[M]}$] by Section \ref{nmb:4.3} and Theorem \ref{nmb:6.1}.
\item[$C^\infty$] by \cite[2.14.3 and 5.26]{KM97}%
 \item[$C^\omega$] by \cite[11.11 and 11.12]{KM97} or by Theorems \ref{nmb:6.1} and \ref{nmb:7.2}, 
\item[$L$] by \cite[2.14.3 and 5.18]{KM97}%
 \item[$\ell^\infty$] by \cite[2.15, 5.24, and 5.25]{KM97} or \cite[3.6.1 and 3.6.6]{FK88}%
 \item[$\operatorname{\mathcal Lip}^k$] by \cite[4.4.2 and 4.4.7]{FK88}
\end{itemize}

It remains to check that $\tilde f$ is of the appropriate class:
\begin{itemize}
\item[\thetag{1}--\thetag{4}] For $\alpha \in \{(M),\{M\}\}$ and $\beta \in \{(N),\{N\},\infty\}$ we have
\begin{align*}
C^\alpha(W_1,C^\beta(W_2,F)) &\cong
L(\lambda^\alpha(W_1),C^\beta(W_2,F))\quad\text{by Lemma \ref{nmb:8.9}}\\
&\cong C^\beta(W_2,L(\lambda^\alpha(W_1),F))\quad\text{by Lemma \ref{nmb:8.8}, 
\cite[3.13.4 and 5.3]{KM97}}\\
&\cong C^\beta(W_2,C^\alpha(W_1,F))\qquad\text{ by Lemma \ref{nmb:8.9}.}
\end{align*}
\item[\thetag{5}] follows from \thetag{2}, \thetag{3}, and Theorem \ref{nmb:7.2}. 
\item[\thetag{6}] is exactly Lemma \ref{nmb:8.8}. 
\item[\thetag{7}] 
follows from
\thetag{6}, using the free convenient vector spaces $\ell^1(X)$ 
over 
the $\ell^\infty$-space $X$, see 
\cite[5.1.24 or 5.2.3]{FK88}, 
satisfying $\ell^\infty(X,F)\cong L(\ell^1(X),F)$.
\item[\thetag{8}]
follows from
\thetag{6}, using the free convenient vector spaces 
$\lambda^k(X)$ over 
the $\operatorname{\mathcal Lip}^k$-space $X$, see \cite[5.1.24 or 5.2.3]{FK88}, 
satisfying $\operatorname{\mathcal Lip}^k(X,F)\cong L(\lambda^k(X),F)$. 
\qed
\end{itemize}
\end{demo}

\section{Manifolds of \texorpdfstring{$C^{{[M]}}$}{C^[M]}-mappings}\label{nmb:9}

\subsection{Hypothesis}
\label{nmb:9.1}
In this section we assume that
$M=(M_k)$ is log-convex and has moderate growth.
In the Beurling case $C^{[M]} = C^{(M)}$ we also require that $C^\omega\subseteq C^{(M)}$, 
equivalently, 
$M_k^{1/k} \to \infty$ or $M_{k+1}/M_k \to \infty$.

For the equivalence of $C^\omega\subseteq C^{(M)}$ and $M_k^{1/k} \to \infty$, see Proposition \ref{nmb:8.1}(4).
Moreover, $M_k^{1/k} \to \infty$ implies $M_{k+1}/M_k \to \infty$, since $M_k^{1/k}$ is increasing, by log-convexity (see Section \ref{nmb:2.1}), 
and thus $M_{k+1}/M_k \ge M_k^{1/k}$. 
Conversely, if $M_{k+1}/M_k \to \infty$ then for each $n \in \mathbb N$ there is $k_n$ so that $M_k/M_{k-1}\ge n$ for all $k\ge k_n$.
It follows that $M_k/M_{k_n-1} \ge n^{k-k_n+1}$ and thus $M_k^{1/k} \to \infty$.
This is needed for the $C^{(M)}$ inverse function theorem (see Sections \ref{nmb:2.1} and \ref{nmb:9.2}).

\subsection{\label{nmb:9.2}Tools for \texorpdfstring{$C^{[M]}$}{C^[M]}-analysis}
We collect here results which are needed below (see also Section \ref{nmb:2.1}):

\begin{enumerate}
\item[\thetag{1}] On open sets in $\mathbb R^n$, $C^{[M]}$-vector fields have $C^{[M]}$-flows, see \cite{Komatsu80} and \cite{Yamanaka91}. 
\item[\thetag{2}] Between Banach spaces, the $C^{[M]}$ implicit function theorem holds. This is essentially due to \cite{Yamanaka89}, but 
in \cite{Yamanaka89} only the Roumieu case is treated and the $C^{\{M\}}$-conditions are global. 
So we shall indicate briefly how to obtain the result we need (cf.\ \cite{RainerSchindl14}): 
\end{enumerate}

\begin{theorem*}
Let $M=(M_k)$ be log-convex.
In the Beurling case $C^{[M]} = C^{(M)}$ we also assume $M_{k+1}/M_k \to \infty$.
Let $E$, $F$ be Banach spaces, $U \subseteq E$, $V \subseteq F$ open, and $f : U \to V$ a $C^\infty$-diffeomorphism.
We have:
\begin{enumerate}
\item[\thetag{3}] Let $K \subseteq U$ be compact. If $f \in C_K^{[M]}(U,F)$ then $f^{-1} \in C_{f(K)}^{[M]}(V,E)$.
\item[\thetag{4}] If $f \in C^{[M]}(U,F)$ then $f^{-1} \in C^{[M]}(V,E)$.
\end{enumerate} 
\end{theorem*}

\begin{demo}{Proof}
By Proposition \ref{nmb:4.2}(5), \thetag{3} implies \thetag{4}. 
The proof of \cite[Thm.\ 2]{Yamanaka89} with small obvious modifications 
provides a proof of \thetag{3} in the Roumieu case (see also \cite[3.4.5]{Schindl09}).

For the Beurling case let $f \in C_K^{(M)}(U,F)$ and
\[
L_k := \frac{1}{k!} \sup_{x \in K} \|f^{(k)}(x)\|_{L^k(E,F)}.
\]
Then $L \lhd M$ and since $M_{k+1}/M_k \to \infty$
there exists a log-convex sequence $N=(N_k)$ satisfying $N_{k+1}/N_k \to \infty$
and such that $L \le N \lhd M$,
by \cite[Lemma 6]{Komatsu79b}.
\begin{comment}
OLD PROOF:
The Beurling case one uses \cite[Lemma 6]{Komatsu79b}:
\emph{Let $M=(M_k)$ be a weight sequence satisfying $M_{k+1}/M_k \to \infty$.
For each sequence $L=(L_k) \subseteq \mathbb{R}_{>0}$ with $L \lhd M$ 
there exists a weight sequence $N=(N_k)$ satisfying $N_{k+1}/N_k \to \infty$
and such that $L \le N \lhd M$.}
\newline
To see this consider $L=(L_k)$ with
\[
L_k := \frac{1}{k!} \sup_{x \in K} \|f^{(k)}(x)\|_{L^k(E,F)}.
\]
Since $f \in C_K^{(M)}(U,F)$, we have $L \lhd M$. So there exists a weight sequence $N=(N_k)$ such that $L \le N \lhd M$. 
\end{comment}
Thus, $f \in C_K^{\{N\}}(U,F)$ and, by the Roumieu case, $f^{-1} \in C_{f(K)}^{\{N\}}(V,E)$. 
Since $N \lhd M$, we have
$f^{-1} \in C_{f(K)}^{(M)}(V,E)$, by Proposition \ref{nmb:8.1}.
\qed\end{demo}

The $C^{[M]}$ implicit function theorem follows in the standard way.

\subsection{\label{nmb:9.3}\texorpdfstring{$C^{[M]}$}{C^[M]}-manifolds} 
A $C^{[M]}$-manifold is a smooth manifold such that all chart changings are
$C^{[M]}$-mappings. They will be considered with the topology induced by 
the $c^\infty$-topology on the charts.
Likewise for $C^{[M]}$-bundles and $C^{[M]}$ Lie groups. 

A mapping between $C^{[M]}$-manifolds is $C^{[M]}$ if and only if
it maps $C^{[M]}$-plots (i.e., $C^{[M]}$-mappings from open sets (or unit balls) of Banach spaces
into the domain manifold) to such.

Note that any finite dimensional (always assumed paracompact) 
$C^\infty$-manifold admits a
$C^\infty$-diffeomorphic real analytic structure thus also a
$C^{[M]}$-structure.

Maybe, any finite dimensional $C^{[M]}$-manifold admits a
$C^{[M]}$-diffeomorphic real analytic structure.
This would follow from:

\subsection*{Conjecture} 
{\it
Let $X$ be a finite dimensional real analytic manifold. Consider the space $C^{[M]}(X,\mathbb R)$ of 
all $C^{[M]}$-functions on $X$, equipped with the (obvious) Whitney $C^{[M]}$-topology. Then $C^\omega(X,\mathbb R)$ is 
dense in $C^{[M]}(X,\mathbb R)$.
}
\medskip

This conjecture is the analogue of \cite[Proposition 8]{Grauert58}.
It was proved in the non-quasianalytic Beurling case $C^{(M)}$ for $X$ open in $\mathbb{R}^n$ by \cite{Langenbruch03}. 

The proofs of the following results are similar to the proofs given in 
\cite[Section 5]{KMRq}, using other analytical tools. For the convenience of the reader, we give 
full proofs here, sometimes with more details.

\subsection{\label{nmb:9.4}Spaces of \texorpdfstring{$C^{[M]}$}{C^[M]}-sections} 
Let $p:E\to B$ be a $C^{[M]}$ vector bundle (possibly infinite dimensional).
The space $C^{[M]}(B\gets E)$ of all $C^{[M]}$-sections is a convenient
vector space with the structure induced by 
\begin{gather*}
C^{[M]}(B\gets E) \to \prod_\alpha C^{[M]}(u_\alpha(U_\alpha),V)
\\
s\mapsto \on{pr}_2\o \psi_\alpha \o s \o u_\alpha^{-1}
\end{gather*}
where $B\supseteq U_\alpha \East{u_\alpha}{} u_\alpha(U_\alpha)\subseteq W$ is a
$C^{[M]}$-atlas for $B$ which we assume to be modeled on a convenient vector
space $W$, and where
$\psi_\alpha :E|_{U_\alpha}\to U_\alpha\times V$ form a vector bundle atlas over
charts $U_\alpha$ of $B$.

\begin{lemma*}
Assume Hypothesis \ref{nmb:9.1}. 
Let $D$ be the open unit ball in a Banach space. A mapping 
$c:D\to C^{[M]}(B\gets E)$ is a $C^{[M]}$-plot if and only if $c^\wedge: D\times B\to E$ is $C^{[M]}$.
\end{lemma*}

\begin{demo}{Proof}
By the description of the structure on $C^{[M]}(B\gets E)$ we may assume by Lemma \ref{nmb:5.1}
that
$B$ is $c^\infty$-open in a convenient vector space $W$ and that $E=B\times V$. 
Then we have $C^{[M]}(B\gets B\times V)\cong C^{[M]}(B,V)$. Thus the statement follows from the
exponential law, i.e., Theorem \ref{nmb:5.2}. 
\qed\end{demo}

Let $U\subseteq E$ be an open neighborhood of $s(B)$ for a section $s$ 
and let $q:F\to B$ be another vector bundle. 
The set $C^{[M]}(B\gets U)$ of all $C^{[M]}$-sections $s':B\to E$ with $s'(B)\subseteq U$ is 
$c^\infty$-open in the convenient vector space $C^{[M]}(B\gets E)$ if $B$ is compact and thus finite 
dimensional, since then it is open in the coarser compact-open topology. 
An immediate consequence of the lemma is the following: If $U\subseteq E$ is an open
neighborhood of $s(B)$ for a section $s$ and if $f:U\to
F$ is a fiber respecting $C^{[M]}$-mapping where 
$F\to B$ is another vector bundle, then $f_*:C^{[M]}(B\gets U)\to
C^{[M]}(B\gets F)$ is $C^{[M]}$ on the open neighborhood $C^{[M]}(B\gets U)$ of $s$ in
$C^{[M]}(B\gets E)$. We have $(d(f_*)(s)v)_x = d(f|_{U\cap E_x})(s(x))(v(x))$.

\begin{theorem}\label{nmb:9.5}
Assume Hypothesis 
\ref{nmb:9.1}. 
Let $A$ and $B$ be finite dimensional $C^{[M]}$-manifolds with $A$ compact and $B$ equipped with a 
$C^{[M]}$ Riemann metric. 
Then the space $C^{[M]}(A,B)$ of all $C^{[M]}$-mappings $A\to B$ is a
$C^{[M]}$-manifold modeled on convenient vector spaces $C^{[M]}(A\gets f^*TB)$ of
$C^{[M]}$-sections of pullback bundles along $f:A\to B$. 
Moreover, a mapping $c:D\to C^{[M]}(A,B)$ is a $C^{[M]}$-plot if and only if 
$c^\wedge :D\times A\to B$ is $C^{[M]}$.
\end{theorem}

If the $C^{[M]}$-structure on $B$ is induced by a real analytic structure, then there exists a real 
analytic Riemann metric which in turn is $C^{[M]}$.

\begin{demo}{Proof}
$C^{[M]}$-vector fields have $C^{[M]}$-flows by Section \ref{nmb:9.2}; %
applying this to the geodesic spray 
we get the $C^{[M]}$ exponential mapping 
$\exp: TB\supseteq U\to B$ of the Riemann metric,
defined on a suitable open neighborhood of the zero section. 
We may assume that $U$ is chosen in such a way that
$(\pi_B,\exp):U\to B\times B$ is a $C^{[M]}$-diffeomorphism onto
an open neighborhood $V$ of the diagonal, 
by the $C^{[M]}$ inverse function theorem, see Section \ref{nmb:9.2}.

For $f\in C^{[M]}(A,B)$ we consider the pullback vector bundle
\[\xymatrix{
A\times TB & A\times_BTB \ar@{_(->}[l] \ar@{=}[r] & f^*TB \ar[r]^{\pi_B^*f} \ar[d]_{f^*\pi_B} & 
TB \ar[d]^{\pi_B} \\
& & A \ar[r]^f & B
}\]
Then the convenient space of sections $C^{[M]}(A\gets f^*TB)$ is canonically isomorphic to the space 
$C^{[M]}(A,TB)_f:= \{h\in C^{[M]}(A,TB):\pi_B\o h=f\}$ 
via $s\mapsto (\pi_B^*f)\o s$ and $(\on{Id}_A,h)\mapsfrom h$. 
Now let 
\begin{gather*} 
U_f :=\{g\in C^{[M]}(A,B):(f(x),\; g(x))\in V
\text{ for all }x\in A\},\\
u_f:U_f\to C^{[M]}(A\gets f^*TB),\\
u_f(g)(x) = (x,\exp_{f(x)}^{-1}(g(x))) = (x,((\pi_B,\exp)^{-1}\o(f,g))(x)).
\end{gather*}
Then $u_f:U_f \to \{s\in C^{[M]}(A\gets f^*TB): s(A)\subseteq f^*U=(\pi_B^*f)^{-1}(U)\}$ is a bijection with 
inverse 
$u_f^{-1}(s) = \exp\o(\pi_B^*f)\o s$, where we view 
$U \to B$ as a fiber bundle. 
The set $u_f(U_f)$ is $c^\infty$-open in
$C^{[M]}(A\gets f^*TB)$ for the topology described above in Section \ref{nmb:9.4}, since $A$ is compact and
the push forward $u_f$ is $C^{[M]}$, since it respects $C^{[M]}$-plots,
by the lemma in Section \ref{nmb:9.4}.

Now we consider the atlas $(U_f,u_f)_{f\in C^{[M]}(A,B)}$ for
$C^{[M]}(A,B)$. Its chart change mappings are given 
for $s\in u_g(U_f\cap U_g)\subseteq C^{[M]}(A\gets g^*TB)$ by 
\begin{align*} 
(u_f\o u_g^{-1})(s) &= (\on{Id}_A,(\pi_B,\exp)^{-1}\o(f,\exp\o(\pi_B^*g)\o s)) \\
&= (\tau_f^{-1}\o\tau_g)_*(s),
\end{align*}
where $\tau_g(x,Y_{g(x)}) := (x,\exp_{g(x)}(Y_{g(x)}))$
is a $C^{[M]}$-diffeomorphism 
$\tau_g:g^*TB \supseteq g^*U \to (g\times \on{Id}_B)^{-1}(V)\subseteq A\times B$
which is fiber respecting over $A$. 
The chart change $u_f\o u_g^{-1} = (\tau_f^{-1}\o \tau_g)_*$ is defined on an open
subset and it is also $C^{[M]}$, since it respects $C^{[M]}$-plots,
by the lemma in Section \ref{nmb:9.4}. 

Finally for the topology on $C^{[M]}(A,B)$ 
we take the identification topology from this atlas (with the 
$c^\infty$-topologies on the modeling spaces $C^{[M]}(A\gets f^*TB)$), 
which is obviously finer than the
compact-open topology and thus Hausdorff.

The equation $u_f\o u_g^{-1} = (\tau_f^{-1}\o \tau_g)_*$ shows that
the $C^{[M]}$-structure does not depend on the choice of the
$C^{[M]}$ Riemannian metric on $B$.

The statement on $C^{[M]}$-plots follows from the lemma in Section \ref{nmb:9.4}. 
\qed\end{demo}

\begin{corollary}\label{nmb:9.6}
Assume Hypothesis 
\ref{nmb:9.1}. 
Let $A_1,A_2$ and $B$ be finite dimensional $C^{[M]}$-manifolds with $A_1$ and
$A_2$ compact. 
Then composition 
\[
C^{[M]}(A_2,B) \times C^{[M]}(A_1,A_2) \to C^{[M]}(A_1,B), \quad (f,g) \mapsto f\o g
\]
is $C^{[M]}$. 
\end{corollary}

\begin{demo}{Proof}
Composition maps $C^{[M]}$-plots to $C^{[M]}$-plots, so it is $C^{[M]}$. 
\qed\end{demo}

\begin{example} \label{nmb:9.7}
The result in Corollary \ref{nmb:9.6} is best possible in the following sense:
If $N=(N_k)$ is another weakly log-convex sequence such that $C^{[N]} \subsetneq C^{[M]}$
(or equivalently, $\inf (N_k/M_k)^{1/k} = 0$ and $\sup (N_k/M_k)^{1/k} < \infty$), 
then composition 
\[
C^{[M]}(S^1,\mathbb R) \times C^{[M]}(S^1,S^1) \to C^{[M]}(S^1,\mathbb R), \quad (f,g) \mapsto f\o g
\]
is {\bf not} $C^{[N]}$ with respect to the canonical real analytic manifold structures.

Namely,
there exists $f\in
C^{[M]}(S^1,\mathbb R)\setminus C^{[N]}(S^1,\mathbb R)$. 
We consider $f$ as a periodic function $\mathbb R\to \mathbb R$.
The universal covering space of $C^{[M]}(S^1,S^1)$ consists of all 
$2\pi\mathbb Z$-equivariant mappings in $C^{[M]}(\mathbb R,\mathbb R)$, 
namely the space of all $g+\on{Id}_{\mathbb R}$ for $2\pi$-periodic 
$g\in C^{[M]}$. Thus $C^{[M]}(S^1,S^1)$ is a real analytic manifold and 
$t\mapsto (x\mapsto x+t)$ induces a real analytic curve $c$ in $C^{[M]}(S^1,S^1)$. 
But $f_*\o c$ is not $C^{(N)}$
(resp.\ $C^{\{N\}}$) since:
\begin{align*}
\frac{(\partial_t^k|_{t=0}(f_*\o c)(t))(x)}{k!\rho^k N_k} =
\frac{\partial_t^k|_{t=0} f(x+t)}{k!\rho^k N_k} = \frac{f^{(k)}(x)}{k!\rho^k N_k}
\end{align*}
which is unbounded in $k$ for $x$ in a suitable compact set and for some (resp.\ all) $\rho>0$, since
$f\notin C^{(N)}$ (resp.\ 
$f\notin C^{\{N\}}$).	
\end{example}

\begin{theorem}\label{nmb:9.8}
Assume Hypothesis 
\ref{nmb:9.1}. 
Let $A$ be a compact (thus finite dimensional) $C^{[M]}$-manifold.
Then the group $\on{Diff}^{[M]}(A)$ of all $C^{[M]}$-diffeomorphisms of $A$ is an
open subset of the $C^{[M]}$-manifold $C^{[M]}(A,A)$. Moreover, it is
a $C^{[M]}$-regular $C^{[M]}$ Lie group: 
Inversion and composition are $C^{[M]}$. Its Lie algebra
consists of all $C^{[M]}$-vector fields on $A$, with the negative of the usual
bracket as Lie bracket. 
The exponential mapping is $C^{[M]}$. It is not surjective onto any neighborhood
of $\on{Id}_A$. 
\end{theorem}

Following \cite{KM97r}, see also \cite[38.4]{KM97}, 
a $C^{[M]}$-Lie group $G$ with Lie algebra $\mathfrak g=T_eG$ 
is called $C^{[M]}$-regular if the following holds:
\begin{itemize}
\item 
For each $C^{[M]}$-curve 
$X\in C^{[M]}(\mathbb R,\mathfrak g)$ there exists a $C^{[M]}$-curve 
$g\in C^{[M]}(\mathbb R,G)$ whose right logarithmic derivative is $X$, i.e.,
\[
\begin{cases} g(0) &= e \\
\partial_t g(t) &= T_e(\mu^{g(t)})X(t) = X(t).g(t)
\end{cases} 
\]
The curve $g$ is uniquely determined by its initial value $g(0)$, if it
exists.
\item
Put $\on{evol}^r_G(X)=g(1)$, where $g$ is the unique solution required above. 
Then $\on{evol}^r_G: C^{[M]}(\mathbb R,\mathfrak g)\to G$ is required to be
$C^{[M]}$ also. 
\end{itemize}

\begin{demo}{Proof}
The group $\on{Diff}^{[M]}(A)$ is $c^\infty$-open in $C^{[M]}(A,A)$, since 
the $C^\infty$-\-diffeo\-mor\-phism group $\on{Diff}(A)$ is $c^\infty$-open in $C^\infty(A,A)$,
by \cite[43.1]{KM97}, and since $\on{Diff}^{[M]}(A)=\on{Diff}(A)\cap C^{[M]}(A,A)$, by Section \ref{nmb:9.2}. 
So $\on{Diff}^{[M]}(A)$ is a $C^{[M]}$-manifold and composition is $C^{[M]}$, by
Theorem \ref{nmb:9.5} and Corollary \ref{nmb:9.6}. To show that inversion is $C^{[M]}$ let $c$ be a
$C^{[M]}$-plot in $\on{Diff}^{[M]}(A)$. By Theorem \ref{nmb:9.5}, the mapping 
$c^\wedge: D\times A\to A$ is
$C^{[M]}$ and $(\on{inv}\o\, c)^\wedge:D\times A\to A$ satisfies the Banach manifold
implicit equation 
$c^\wedge(t,(\on{inv}\o\, c)^\wedge(t,x))=x$ for $x\in A$. By the Banach
$C^{[M]}$ implicit function theorem, see Section \ref{nmb:9.2},  %
the mapping 
$(\on{inv}\o\, c)^\wedge$ is locally $C^{[M]}$ and thus $C^{[M]}$. By Theorem \ref{nmb:9.5}
again, $\on{inv}\o\, c$ is a $C^{[M]}$-plot in $\on{Diff}^{[M]}(A)$.
So $\on{inv}:\on{Diff}^{[M]}(A)\to \on{Diff}^{[M]}(A)$ is $C^{[M]}$. 
The Lie algebra of $\on{Diff}^{[M]}(A)$ is the convenient vector space of all 
$C^{[M]}$-vector fields on $A$, with the negative of the usual Lie bracket
(compare with the proof of \cite[43.1]{KM97}). 

To show that $\on{Diff}^{[M]}(A)$ is a $C^{[M]}$-regular Lie group, we choose a
$C^{[M]}$-plot in the space of $C^{[M]}$-curves in the Lie algebra of all $C^{[M]}$
vector fields on $A$, that is $c:D\to C^{[M]}(\mathbb R,C^{[M]}(A\gets TA))$. By the lemma in
Section \ref{nmb:9.4}, the plot $c$ corresponds to a $(D\times \mathbb R)$-time-dependent $C^{[M]}$
vector field $c^{\wedge \wedge }:D\times\mathbb R\times A\to TA$. Since $C^{[M]}$-vector
fields have $C^{[M]}$-flows and since $A$ is compact, 
$\on{evol}^r(c^\wedge (s))(t) = \on{Fl}^{c^\wedge (s)}_t$ is $C^{[M]}$ in all
variables, by Section \ref{nmb:9.2}. %
Thus $\on{Diff}^{[M]}(A)$ is a $C^{[M]}$-regular $C^{[M]}$ Lie group.

The exponential mapping is $\on{evol}^r$ applied to constant curves in the
Lie algebra, i.e., it consists of flows of autonomous $C^{[M]}$ vector fields. 
That the exponential mapping is not surjective onto any $C^{[M]}$-neighborhood of
the identity follows from \cite[43.5]{KM97} for $A=S^1$. This example can
be embedded into any compact manifold, see \cite{Grabowski88}. 
\qed\end{demo}

%\bibliography{../../references/biblio}
%\bibliographystyle{amsplain}

\def\cprime{$'$}
\providecommand{\bysame}{\leavevmode\hbox to3em{\hrulefill}\thinspace}
\providecommand{\MR}{\relax\ifhmode\unskip\space\fi MR }
\providecommand{\MRhref}[2]{%
  \href{http://www.ams.org/mathscinet-getitem?mr=#1}{#2}
}
\providecommand{\href}[2]{#2}

\end{document}